\documentclass[11pt]{article}
\usepackage[a4paper,margin=1in]{geometry}
\usepackage{amsmath,amssymb,amsthm,mathtools,bm}
\usepackage{microtype}
\usepackage{enumitem}
\usepackage[colorlinks=true,linkcolor=blue,citecolor=blue,urlcolor=blue]{hyperref}
\usepackage[nameinlink,capitalize]{cleveref}

\newtheorem{theorem}{Theorem}[section]
\newtheorem{proposition}[theorem]{Proposition}
\newtheorem{lemma}[theorem]{Lemma}
\newtheorem{corollary}[theorem]{Corollary}
\theoremstyle{remark}
\newtheorem{remark}[theorem]{Remark}

\crefname{proposition}{Proposition}{Propositions}
\crefname{lemma}{Lemma}{Lemmas}
\crefname{corollary}{Corollary}{Corollaries}

\newcommand{\R}{\mathbb R}
\newcommand{\C}{\mathbb C}
\newcommand{\N}{\mathbb N}
\newcommand{\Z}{\mathbb Z}
\newcommand{\T}{\mathbb T}

\newcommand{\Om}{\Omega}
\newcommand{\eps}{\varepsilon}

\newcommand{\dd}{\,\mathrm d}
\newcommand{\supp}{\operatorname{supp}}
\newcommand{\Lip}{\mathrm{Lip}}
\newcommand{\Id}{\mathrm{Id}}

\newcommand{\diag}{\operatorname{diag}}

\newcommand{\dist}{\operatorname{dist}}

\title{Almost-sure quenched KAM tori for nonlinear wave equations
with a spatial white-noise potential}
\author{
Yingdu Dong\thanks{Affiliation: School of Mathematical Sciences, Fudan University, Shanghai, China. Email: 202131130009@mail.bnu.edu.cn.}
\and
Wenwen Jian\thanks{Affiliation: School of Mathematics, Physics and Statistics,
Shanghai Polytechnic University, Shanghai, China. Email: wwjian@sspu.edu.cn. Funding: NNSFC(11201392).}
\and
Xiaoping Yuan\thanks{Affiliation: School of Mathematical Sciences, Fudan University, Shanghai, China. Email: xpyuan@fudan.edu.cn. Funding: NNSFC(12371189).}
}
\date{}
\newcommand{\cE}{\mathcal E}
\newcommand{\cL}{\mathcal L}

\newcommand{\cF}{\mathcal F}
\newcommand{\Hyp}{\mathcal H_-}
\newcommand{\Ell}{\mathcal E_+}

\begin{document}
\maketitle

\begin{abstract}
Let
\[
 A_\omega=-\partial_x^2+\rho\,\partial_xB(x,\omega),
 \qquad \rho\ne0,
\]
be the Dirichlet Schr\"odinger operator on $(0,\pi)$ with a spatial
Gaussian white-noise potential, realized pathwise via quasi-derivatives.
Here $B$ is Brownian motion in the spatial variable $x$, and
$\partial_xB$ denotes its distributional derivative.  We consider the
nonlinear wave equation
\[
 u_{tt}+A_\omega u+f(u)=0,
\]
where $f$ is odd and real analytic near zero, with
$f(u)=\kappa u^3+O(u^5)$ and $\kappa\ne0$.
We show that, on a single event of probability one, no nontrivial
finitely supported integer combination of the positive wave frequencies
is an integer, and the quartic twist matrix is nonsingular for every
finite tangential set of positive modes.  After fixing a path in this
event, we construct a partial Birkhoff normal form of the Hamiltonian
up to degree ten.
For every finite nonempty tangential set $J$ of cardinality $b$
associated with positive eigenvalues, we obtain a Cantor family of real-analytic,
small-amplitude invariant $b$-tori in
$H^1_0(0,\pi;\R)\times L^2(0,\pi;\R)$.  The family is parametrized
by a Cantor subset of $[\nu,2\nu]^b$ whose relative measure tends to
one as $\nu\to0$.  The tori are normally elliptic when $A_\omega$ is
positive.  Otherwise, they have $h$ stable and $h$ unstable normal
directions, where $h$ is the number of negative eigenvalues of $A_\omega$.
\end{abstract}

\medskip
\noindent\textbf{Keywords.}
Nonlinear wave equation; spatial white noise; quenched KAM theory;
partial Birkhoff normal form; Gaussian zero sets; hyperbolic--elliptic tori.

\medskip
\noindent\textbf{2020 Mathematics Subject Classification.}
37K55, 35L71, 60H40, 34L40.

\tableofcontents

\section{Introduction}
We consider the real nonlinear wave equation
\begin{equation}\label{eq:NLW}
 u_{tt}+A_\omega u+f(u)=0,\qquad
 A_\omega=-\partial_x^2+\rho\,\partial_xB(x,\omega),\qquad 0<x<\pi,
\end{equation}
with Dirichlet boundary conditions, where $\rho\in\R\setminus\{0\}$
and $B=(B_x)_{0\le x\le\pi}$ is standard real Brownian motion.
The derivative $\partial_xB$ is understood in the distributional sense.
We assume that $f$ is odd and real analytic in a neighborhood of zero,
with
\begin{equation}\label{eq:nonlinearity}
 f(s)=\kappa s^3+O(s^5),\qquad \kappa\in\R\setminus\{0\}.
\end{equation}

This equation describes wave propagation in a
spatially disordered one-dimensional medium, where the real potential
represents spatial fluctuations of the medium.  Gaussian white noise
provides a model for fluctuations whose correlation length is small
compared with the observation scale.  Its covariance is
\[
 \mathbb E\bigl[\langle\rho\partial_xB,\varphi\rangle
                    \langle\rho\partial_xB,\psi\rangle\bigr]
       =\rho^2\int_0^\pi \varphi(x)\psi(x)\dd x,
 \qquad \varphi,\psi\in C_c^\infty(0,\pi;\R).
\]
For further discussion of the physical background of the white-noise
potential, we refer to our previous
work~\cite{DongJianYuanSpectral2026}.

Schr\"odinger operators with white-noise potentials were studied by
Fukushima and Nakao~\cite{FukushimaNakao1977}; more recent work includes
the construction and localization theory on the
line~\cite{DumazLabbe2024}.
On the finite interval considered here, the primitive $\rho B$ is
almost surely continuous.  The quasi-derivative theory of singular
Sturm--Liouville operators therefore gives a pathwise self-adjoint
realization with compact resolvent and form domain $H^1_0(0,\pi)$;
see~\cite{SavchukShkalikov2003,HrynivMykytyuk2006,HrynivMykytyuk2012}.
The finite-interval spectral estimates used below were established
in~\cite{DongJianYuanSpectral2026}.

Infinite-dimensional KAM theory has been used to construct
quasi-periodic solutions of nonlinear wave equations with both
parameter-dependent and prescribed potentials.  Wayne
\cite[Sections~4 and~9]{Wayne1990} introduced a Gaussian measure on
a class of regular potentials through inverse spectral coordinates.
The spectral coordinates serve as external parameters, and invariant
tori are constructed after excluding a set of potentials of small
measure.  For equations with a fixed linear part, the tangential
actions can instead provide the parameters required by KAM theory.
P\"oschel \cite{PoschelWave1996} used this approach for the fixed-mass
wave equation.  Yuan extended the KAM theorem to finitely many
hyperbolic normal modes \cite{YuanHyperbolic2003} and treated
prescribed smooth potentials with nonzero mean
\cite{YuanWave2006}, with tangential modes chosen sufficiently high
and suitably separated.

Random potentials also arise in KAM theory for Hamiltonian lattices;
see Yuan's survey \cite{YuanSurvey2014} for the role of the linear
spectral properties in systems with short-range interactions.
More recently, Geng, Sun, and Wang \cite{GengSunWangWave2026}
constructed quasi-periodic solutions of a nonlinear random wave
equation on $\Z^d$.  Their KAM theorem in
\cite{GengSunWangKAM2026} treats discrete nonlinear Schr\"odinger
and wave equations with spatially decaying coupling and also yields
linear stability.  In these results, a random potential is first
selected from a suitable set and then fixed, after which the small
amplitudes serve as internal parameters.  Our construction follows
the same order: probabilistic selection of the potential precedes
the KAM construction with action parameters.

The main difference concerns the probability law under which this
selection is made.  In Wayne's formulation, the measure is introduced
through spectral coordinates; in the lattice setting, randomness is
specified by the site potentials.  Here the probability law is fixed
by the spatial white-noise model itself.  The sample space is the
Wiener space of the Brownian primitive, and the spectral data are
determined by the corresponding path.  In particular, the
eigenvalues and eigenfunctions are correlated nonlinear functionals
of that path, not independently prescribed spectral parameters.
Thus the nonresonance and twist conditions must be established under
the given Wiener measure.  For this purpose, we use the analytic
dependence of the spectral data and the zero-set theorem on Wiener
space proved in our work on the cubic nonlinear Schr\"odinger
equation \cite{DongJianYuanNLS2026}, together with the almost-sure
quartic twist condition obtained there.  We prove that the required
wave-frequency nonresonance and twist conditions hold on one event
of probability one for every finite tangential set of positive
modes.  Once a path in this event is fixed, we construct a higher-order partial Birkhoff
normal form and use the tangential actions as parameters in
P\"oschel's theorem \cite{Poschel1996} and Yuan's extension
\cite{YuanHyperbolic2003}.

The operator $A_\omega$ may have finitely many negative eigenvalues.
They give hyperbolic directions in the wave equation.  The tangential
modes are therefore chosen from the positive spectrum, while the
negative modes are retained as normal variables. Write
\[
 \Ell(\omega)=\{n:\lambda_n(\omega)>0\},\qquad
 \Hyp(\omega)=\{n:\lambda_n(\omega)<0\},\qquad
 h(\omega)=|\Hyp(\omega)|.
\]
The set $\Hyp$ is finite.  For $n\in\Ell$ put
$\mu_n=\sqrt{\lambda_n}$, and choose real $L^2$-normalized eigenfunctions
$\phi_n$.  Given $J=\{j_1,\ldots,j_b\}\subset\Ell$, the linear rotational
torus with actions $\xi\in(0,\infty)^b$ is parametrized by
\begin{equation}\label{eq:linear-torus}
 \begin{split}
 u_J(\theta)&=\sum_{a=1}^b\sqrt{2\xi_a/\mu_{j_a}}
                  \cos\theta_a\,\phi_{j_a},\\
 v_J(\theta)&=-\sum_{a=1}^b\sqrt{2\mu_{j_a}\xi_a}
                  \sin\theta_a\,\phi_{j_a}.
 \end{split}
\end{equation}

\begin{theorem}[Almost-sure quenched KAM tori]\label{thm:main}
Fix $\rho\in\R\setminus\{0\}$ and an odd nonlinearity $f$, real analytic
near zero, satisfying \eqref{eq:nonlinearity}.  There is a measurable
event $\Omega_*$ of probability one, independent of $J$ and $f$, on
which $0\notin\operatorname{spec}(A_\omega)$
and the following assertions hold.
For every $\omega\in\Omega_*$, every finite nonempty
$J\subset\Ell(\omega)$, and every $0<\beta<1/6$, there exist
$\nu_0=\nu_0(\omega,f,J,\beta)>0$ and $C=C(\omega,f,J,\beta)$ such that,
for $0<\nu<\nu_0$, a Cantor set $\Pi_\nu^*\subset\Pi_\nu=[\nu,2\nu]^b$
satisfies
\begin{equation}\label{eq:relative-measure}
 \frac{|\Pi_\nu\setminus\Pi_\nu^*|}{|\Pi_\nu|}\le C\nu^\beta.
\end{equation}
For each $\xi\in\Pi_\nu^*$, \eqref{eq:NLW} has a real-analytic invariant
$b$-torus in $H^1_0(0,\pi;\R)\times L^2(0,\pi;\R)$ carrying
quasi-periodic motion.  Its embedding is an $o_{\omega,f,J,\beta}(\sqrt\nu)$
deformation of \eqref{eq:linear-torus} in the energy norm.
The normal linearized system is reducible to constant elliptic
oscillators and $h(\omega)$ hyperbolic pairs.  In particular, if
$\lambda_1(\omega)>0$, the tori are normally elliptic and linearly stable
in the normal directions.  If $h(\omega)>0$, they have $h(\omega)$
stable and $h(\omega)$ unstable normal directions.

The trajectories satisfy
\[
 u\in C(\R;H^1_0)\cap C^1(\R;L^2)\cap C^2(\R;H^{-1})
\]
and solve \eqref{eq:NLW} in $H^{-1}$, equivalently as a mild energy
solution.  In addition,
\[
 (u,u_t)\in C^\omega\bigl(\R;H^1_0(0,\pi;\R)\times L^2(0,\pi;\R)\bigr).
\]
Each trajectory extends holomorphically, with values in the complexified
energy space, to a strip of positive width about the real time axis.
\end{theorem}

\begin{remark}[Choice of tangential modes]
The event $\Omega_*$ is selected before the tangential set.  Hence every
finite set of positive modes is admissible on the same event, with no
high-mode or separation condition.  The almost-sure nonresonance and twist conditions permit this choice.
After the path is fixed, the tangential actions serve as parameters;
the amplitude threshold and the constants in the estimates depend on
the realization and on $J$.
\end{remark}

The proof has two main parts.  The first concerns wave-frequency
nonresonance.  Nonresonance of the eigenvalues $\lambda_n$ does not
imply nonresonance of their positive square roots
$\mu_n=\sqrt{\lambda_n}$.  Moreover, the differences $\mu_n-\mu_m$
approach the integers $n-m$ at high energies.  Algebraic elimination of
the square roots reduces the relevant relations to zero sets of
analytic functions on Wiener space.  This excludes, on one event of
probability one, every nontrivial relation between finitely many positive
wave frequencies and an integer.

The second part is the normal-form construction for the fixed operator.
The nonlinear Hamiltonian has a one-smoothing vector field in the
spectral energy coordinates.  We prove that the solutions of the
homological equations have the same regularity.  For quadratic terms
with elliptic difference divisors, we first solve the equation with
$\mu_n$ replaced by $n$ by Fourier averaging, and then include
$\mu_n-n$ by a convergent Neumann series.  Sum divisors are treated
by a Laplace integral.  At each step we normalize the terms of normal
degree at most two and retain the higher normal degrees in the perturbation.

The normalization order is determined by the measure estimate used in
the KAM theorem.  The asymptotics
$\mu_n=n+O_{\omega,\eta}(n^{-\eta})$, $0<\eta<1$, give the exponent
$\eta/(1+\eta)<1/2$ in that estimate.  With the perturbation bounds
proved below, normalization at degrees $4,6,8,10$ makes the KAM smallness
condition compatible with a relative measure tending to one.
The higher resonant action terms are included in the integrable
Hamiltonian, so its frequencies depend polynomially on the actions.
We verify the frequency hypotheses and the uniformity of the measure
constants on shrinking action boxes.

The rest of this paper is organized as follows. Section~2 recalls the operator and spectral results from
\cite{DongJianYuanNLS2026}.  Section~3 establishes wave-frequency
nonresonance and identifies the quartic twist matrix.  Section~4 proves the
smoothing estimates and constructs the partial normal form.
Section~5 verifies the KAM hypotheses and proves \cref{thm:main}.
Appendix~\ref{app:KAM}
states the KAM results and measure estimates used in Section~5.

\section{Preliminaries}
In this section, we collect the operator-theoretic and spectral results used in the proof of the main theorem.
\subsection{The operator and the energy space}
For a real $\sigma\in L^2(0,\pi)$, define
\[
 y^{[1]}=y'-\sigma y,\qquad
 \ell_\sigma y=-(y^{[1]})'-\sigma y^{[1]}-\sigma^2y.
\]
The Dirichlet realization $A_\sigma$ has domain
\[
 \{y\in L^2:y,y^{[1]}\in AC[0,\pi],\quad
       \ell_\sigma y\in L^2,\quad y(0)=y(\pi)=0\}.
\]
Its closed form on $H^1_0(0,\pi)$ is
\begin{equation}\label{eq:form}
 \mathfrak a_\sigma[u,v]=\int_0^\pi u'\overline{v'}\dd x
 -\int_0^\pi\sigma(u'\overline v+u\overline{v'})\dd x.
\end{equation}
This is the distributional interpretation of $-\partial_x^2+\sigma'$.
The standard distributional Sturm--Liouville theory gives
self-adjointness, lower semiboundedness, compact resolvent, and simple
Dirichlet eigenvalues \cite{SavchukShkalikov2003,HrynivMykytyuk2006,HrynivMykytyuk2012}.

The form estimates in \cite[Section~2.1]{DongJianYuanNLS2026} give, for
some sufficiently large $c_\sigma>0$,
\begin{equation}\label{eq:form-norm}
 \mathfrak a_\sigma[u,u]+c_\sigma\|u\|_2^2\asymp_\sigma\|u\|_{H^1}^2.
\end{equation}

For $\sigma=\rho B(\omega)$, write $A_\omega=A_\sigma$ and
\[
 A_\omega\phi_n=\lambda_n\phi_n,\qquad
 \lambda_1<\lambda_2<\cdots\longrightarrow+\infty,\qquad
 \|\phi_n\|_2=1.
\]
The spectral theorem applied to the shifted form gives
\begin{equation}\label{eq:spectral-norm}
 \left\|\sum_n u_n\phi_n\right\|_{H^1}^2
 \asymp_\omega\sum_n(c_\omega+\lambda_n)|u_n|^2.
\end{equation}
For each fixed path, choose the signs of the real eigenfunctions
arbitrarily and keep that choice throughout the construction.

\subsection{High-energy wave frequencies}
\begin{proposition}\label{prop:spectral}
On an event $\Omega_{\rm sp}$ of probability one, for all $\epsilon>0$,
\begin{equation}\label{eq:spectral-estimates}
 \lambda_n=n^2+O_{\omega,\eps}(n^\eps),\qquad
 \sup_n\|\phi_n\|_\infty<\infty.
\end{equation}
For every $0<\eta<1$, on the positive spectral tail,
\begin{equation}\label{eq:wave-asymptotic}
 \mu_n=n+O_{\omega,\eta}(n^{-\eta}),\qquad
 \frac{\mu_n-\mu_m}{n-m}=1+O_{\omega,\eta}(m^{-\eta})\quad(n>m),
\end{equation}
uniformly for sufficiently large $m$.
\end{proposition}
\begin{proof}
The estimates in \eqref{eq:spectral-estimates} hold on the common event
constructed in \cite[Section~2]{DongJianYuanNLS2026}, using
\cite{DongJianYuanSpectral2026}.

Since $\lambda_n>0$ eventually,
\[
 \mu_n-n=\frac{\lambda_n-n^2}{\sqrt{\lambda_n}+n}
       =O(n^{-\eta})
\]
by taking $\eps=1-\eta$.  If $r_n=\mu_n-n$, then, for $n>m$,
\[
 \left|\frac{r_n-r_m}{n-m}\right|
 \le C\frac{n^{-\eta}+m^{-\eta}}{n-m}\le 2Cm^{-\eta}.
\]
This proves the uniform gap estimate.  In particular consecutive
positive wave-frequency gaps tend to $1$.
\end{proof}

\subsection{Analytic spectral functionals}
Let $E=C_0([0,\pi];\R) =\{B\in C([0,\pi];\R):B(0)=0\}$ with the uniform norm, and let $\mathbb W$ be
classical Wiener measure on $E$.  For deterministic $B\in E$ write
$A(B)=A_{\rho B}$.  By \cite[Section~2]{DongJianYuanNLS2026},
the ordered eigenvalues are continuous and locally real analytic on
$E$, and locally one can choose normalized real eigenfunctions analytic
with values in $H^1_0$.  In particular,
\begin{equation}\label{eq:Cmn}
 C_{mn}(B)=\int_0^\pi\phi_m(B,x)^2\phi_n(B,x)^2\dd x
\end{equation}
is independent of the eigenfunction signs, globally well defined, and
locally real analytic.  The positive square root $\mu_n$ is locally real analytic on
$\{\lambda_n>0\}$.  To study relations among these square roots, we
eliminate them algebraically and apply the zero-set theorem to functions
defined on all of $E$.

\section{Almost-sure nonresonance and twist}\label{sec:probability}
We use the following consequence of
\cite[Theorem~3.1]{DongJianYuanNLS2026}: a locally real-analytic
function $f:E\to\R$ which is not identically zero has a Wiener-null
zero set.  Analyticity is understood in the Banach-space sense, with
locally norm-convergent series of continuous homogeneous polynomials.
The functions to which we apply this theorem are defined on the whole
connected space $E$.

\subsection{Nonresonance of the wave frequencies}
\begin{proposition}\label{prop:radicals}
There is a measurable event $\Omega_{\rm ar}$ of probability one such
that $\lambda_n\ne0$ for every $n$, and
\begin{equation}\label{eq:affine-nr}
 a+\sum_{n\in S}c_n\sqrt{\lambda_n}\ne0
\end{equation}
for every finite $S\subset\Ell$, every $c\in\Z^S$, and every $a\in\Z$
with $(a,c)\ne0$.
\end{proposition}
\begin{proof}
For $c\ne0$ we may take $S=\supp c$.  Introduce indeterminates
$t=(t_n)_{n\in S}$ and form the polynomial
\[
 \mathcal Q_{a,c}(t)=
 \prod_{\epsilon\in\{\pm1\}^S}
       \left(a+\sum_{n\in S}\epsilon_n c_n t_n\right).
\]
Changing the sign of a single $t_n$ permutes the factors.  Consequently
only even powers of every $t_n$ occur, and there is a real polynomial
$Q_{a,c}$ such that $\mathcal Q_{a,c}(t)=Q_{a,c}((t_n^2)_{n\in S})$.
The function
\[
 F_{a,c}(B)=Q_{a,c}((\lambda_n(B))_{n\in S})
\]
is therefore globally defined and locally real analytic on $E$,
including where some eigenvalues are negative.

We prove that it is not identically zero.  On the deterministic path
$B_s(x)=sx/\rho$, $s>0$, the potential is the constant $s$ and
$\lambda_n(B_s)=n^2+s$.  Every factor of $F_{a,c}(B_s)$ has the form
\[
 a+\sum_{n\in S}d_n\sqrt{s+n^2},\qquad d_n=\epsilon_nc_n.
\]
Suppose this analytic function vanished for all $s>0$.  Division by
$\sqrt{s}$ and passage to the limit give $\sum d_n=0$.
The remaining square-root sum is then $O(s^{-1/2})$, so $a=0$.
For $s>\max_{n\in S}n^2$, the convergent expansion
\[
 \sum_{n\in S}d_n\sqrt{s+n^2}
 =\sqrt{s}\sum_{r\ge0}\binom{1/2}{r}s^{-r}
                         \sum_{n\in S}d_n n^{2r}
\]
has all coefficients zero.  In particular,
\[
 \sum_{n\in S}d_n n^{2r}=0,\qquad r=0,1,\ldots,|S|-1.
\]
All binomial coefficients involved are nonzero.  The Vandermonde matrix
on the distinct numbers $n^2$ is invertible, forcing every $d_n=0$,
a contradiction.  Thus each factor is a nontrivial real-analytic function of $s>0$ and
has a discrete zero set.  Their finite product is therefore nonzero for
some $s>0$.

By \cite[Theorem~3.1]{DongJianYuanNLS2026}, $F_{a,c}^{-1}(0)$ is Wiener null.
Whenever all $\lambda_n$, $n\in S$, are positive, a relation
$a+\sum c_n\sqrt{\lambda_n}=0$ would make the all-positive factor of
$\mathcal Q_{a,c}$ vanish.  Removing the null sets for the countably
many pairs $(a,c)$, with $c$ finitely supported, proves
\eqref{eq:affine-nr}.  If $c=0$, the assertion
reduces to $a\ne0$.  Intersecting with the nonresonance event of
\cite[Theorem~3.3]{DongJianYuanNLS2026} also excludes zero eigenvalues.
Taking the inverse image under the Brownian-path random variable gives
the claimed event on the original sample space.
\end{proof}

For each nonzero integer coefficient vector $(c_j)_{j\in J}$,
\eqref{eq:affine-nr} gives
$C=\sum_{j\in J}c_j\mu_j\notin\Z$ and hence
$\dist(C,\Z)>0$.  This separation from integers is used below for
the divisors $C+\mu_n-\mu_m$, since
$\mu_n-\mu_m=(n-m)+o(1)$ when both indices tend to infinity.
Only finitely many tangential combinations occur in the finite
normal-form construction, so their distances from $\Z$ have a positive
path-dependent minimum.

\subsection{The quartic twist matrix}
For $J=\{j_1,\ldots,j_b\}\subset\N$, set
\[
 (G_J)_{aa}=C_{j_aj_a},\qquad
 (G_J)_{ac}=2C_{j_aj_c}\quad(a\ne c).
\]
By \cite[Theorem~3.4]{DongJianYuanNLS2026}, there is one measurable
full-probability event $\Omega_{\rm tw}$ on which $G_J$ is invertible
for every finite nonempty $J$.  For $J\subset\Ell$, the Hessian $A_J$ of the quartic tangential
action polynomial is determined by the cubic coefficient $\kappa$
in \eqref{eq:nonlinearity}:
\begin{equation}\label{eq:twist}
 A_J=\frac{3\kappa}{4}D_J^{-1}G_JD_J^{-1},\qquad
 D_J=\diag(\mu_{j_1},\ldots,\mu_{j_b}).
\end{equation}
The coefficient calculation is given in \cref{sec:coefficients}.
Since $D_J$ is invertible, $A_J$ is nonsingular on this event.

We now fix a realization
\begin{equation}\label{eq:universal-event}
 \omega\in\Omega_*:=\Omega_{\rm sp}\cap\Omega_{\rm ar}\cap\Omega_{\rm tw}.
\end{equation}
The dependence of spectral quantities on $\omega$ is henceforth omitted
from the notation.  Constants may depend on this fixed path, on the
nonlinearity $f$, and on $J$.  The event $\Omega_*$ depends only on the
random operator and is independent of $f$ and the subsequent choice
of $J\subset\Ell$.

\section{Partial Birkhoff normal form}\label{sec:BNF}
\subsection{Hamiltonian formulation and smoothing estimates}
Let $R_f>0$ be a convergence radius for the Taylor series of $f$ at
zero, and define its even primitive by
\[
 \cF(s)=\int_0^s f(\tau)\dd\tau
       =\frac\kappa4 s^4+O(s^6).
\]
The real Hamiltonian of \eqref{eq:NLW} is
\begin{equation}\label{eq:H-physical}
 H(u,v)=\frac12\|v\|_2^2+\frac12\mathfrak a_{\rho B}[u,u]
                       +\int_0^\pi \cF(u(x))\dd x,
 \qquad v=u_t.
\end{equation}
It is real analytic on a neighborhood of zero in $H^1_0\times L^2$,
chosen so that $\|u\|_\infty<R_f$.
Write $u_n=\langle u,\phi_n\rangle$ and $v_n=\langle v,\phi_n\rangle$.
For $n\in\Ell$, introduce the complex coordinates $q_n$ and
$\widetilde q_n=\overline{q_n}$ by
\begin{equation}\label{eq:ell-coordinates}
 q_n=\frac{\sqrt{\mu_n}u_n+i v_n/\sqrt{\mu_n}}{\sqrt2},\qquad
 u_n=\frac{q_n+\widetilde q_n}{\sqrt{2\mu_n}},\qquad
 v_n=-i\sqrt{\frac{\mu_n}{2}}(q_n-\widetilde q_n).
\end{equation}
For $h\in\Hyp$, put $\chi_h=\sqrt{-\lambda_h}>0$ and
\begin{equation}\label{eq:hyp-coordinates}
 a_h=\frac{v_h/\sqrt{\chi_h}+\sqrt{\chi_h}u_h}{\sqrt2},\qquad
 b_h=\frac{v_h/\sqrt{\chi_h}-\sqrt{\chi_h}u_h}{\sqrt2}.
\end{equation}
Thus $u_h=(a_h-b_h)/\sqrt{2\chi_h}$ and
$v_h=\sqrt{\chi_h/2}(a_h+b_h)$.  The symplectic form becomes
\begin{equation}\label{eq:symp}
 \Om=i\sum_{n\in\Ell}\dd q_n\wedge\dd\widetilde q_n
                      +\sum_{h\in\Hyp}\dd a_h\wedge\dd b_h.
\end{equation}
We use $\iota_{X_F}\Om=\dd F$ and $\{F,G\}=\dd F[X_G]$.
In particular, $X_F^{q}=-i\partial_{\widetilde q}F$,
$X_F^{\widetilde q}=i\partial_qF$, $X_F^a=\partial_bF$, and
$X_F^b=-\partial_aF$.  The Hamiltonian is $H=H_2+P$, with
\begin{equation}\label{eq:H2}
 H_2=\sum_{n\in\Ell}\mu_nq_n\widetilde q_n+
                       \sum_{h\in\Hyp}\chi_h a_hb_h.
\end{equation}
Writing $f(s)=\sum_{r\ge2}a_{2r-1}s^{2r-1}$, with $a_3=\kappa$,
the nonlinear Hamiltonian has the expansion
\begin{equation}\label{eq:nonlinear-Hamiltonian}
 P(w)=\int_0^\pi\cF(u(w,x))\dd x
     =\sum_{r\ge2}P_{2r}(w),\qquad
 P_{2r}(w)=\frac{a_{2r-1}}{2r}\int_0^\pi u(w,x)^{2r}\dd x.
\end{equation}
For analytic estimates, we regard $q,\widetilde q,a,b$ as independent
complex variables.  Real-valued wave data correspond to
$\widetilde q_n=\overline{q_n}$ and $a_h,b_h\in\R$.

For $p\ge0$, let
\[
 h^p(\Ell)=\left\{q:\|q\|_p^2=
                 \sum_{n\in\Ell}\langle n\rangle^{2p}|q_n|^2<\infty\right\},
 \quad \langle n\rangle=(1+n^2)^{1/2},
\]
and let $\cE^p=h^p(\Ell)\times h^p(\Ell)\times\C^{2h(\omega)}$,
with the sum norm.  The energy space $H^1_0(0,\pi;\R)\times L^2(0,\pi;\R)$
corresponds to the points of
$\cE^{1/2}$ satisfying these reality conditions.  Indeed, \eqref{eq:spectral-norm} and
$\mu_n\asymp n$ give the norm equivalence on positive modes, and the
remaining transformations are finite dimensional.  The quadratic functional $H_2$ is
continuous on $\cE^{1/2}$, although $X_{H_2}$ is unbounded there.
All brackets with $H_2$ below mean $\dd H_2[X_F]$ for a bounded $X_F$.

For later use, denote the spectral reconstruction map by
$\mathcal Uf=\sum_n f_n\phi_n$ and let $L=\diag(\mu_n)_{n\in\Ell}$.
The positive-mode displacement map is
\[
 \mathcal S(q,\widetilde q)
       =2^{-1/2}\mathcal U L^{-1/2}(q+\widetilde q).
\]
Here $L^{-1/2}:h^{1/2}(\Ell)\to h^1(\Ell)$ and
$\mathcal U:h^1\to H^1_0$ are bounded.  Conversely,
$L^{-1/2}\mathcal U^{-1}$ maps $H^1_0$ to $h^{3/2}(\Ell)$,
where restriction to the positive modes is understood.
The displacement map contains one factor $L^{-1/2}$, and
differentiation of the Hamiltonian contributes another.  Together with
the product estimate in $H^1_0$, these factors give the following
smoothing estimate.  We call a vector field one-smoothing if it is
analytic from a neighborhood of zero in $\cE^{1/2}$ to $\cE^{3/2}$.

\begin{lemma}[Smoothing of the nonlinear vector field]\label{lem:smoothing}
The vector field $X_P$ is real analytic from a neighborhood of zero in
$\cE^{1/2}$ to $\cE^{3/2}$ and satisfies
\begin{equation}\label{eq:nonlinear-smooth}
 \|X_P(w)\|_{3/2}\le C\|w\|_{1/2}^3,\qquad
 \|DX_P(w)\|_{\cL(\cE^{1/2},\cE^{3/2})}
                         \le C\|w\|_{1/2}^2.
\end{equation}
Moreover, the series in \eqref{eq:nonlinear-Hamiltonian} converges
locally uniformly on a sufficiently small complex energy ball, as does
the series of its Hamiltonian vector fields with values in
$\cE^{3/2}$.  Each $X_{P_{2r}}$ is a continuous homogeneous polynomial
of degree $2r-1$.
\end{lemma}
\begin{proof}
The linear map $w\mapsto u$ defined by \eqref{eq:ell-coordinates} and
\eqref{eq:hyp-coordinates} is bounded from $\cE^{1/2}$ to $H^1_0$.
The Banach algebra estimate on $H^1(0,\pi)$ gives, for $j\ge1$,
\[
 \|u^j\|_{H^1}\le C_{\rm alg}^{j-1}\|u\|_{H^1}^j.
\]
For sufficiently small $\|u\|_{H^1}$, the Taylor series of $f(u)$ and
its derivative therefore converge in $H^1$.  Since $f(0)=0$, the
composition preserves the zero boundary traces.  It defines a
holomorphic map on a complex neighborhood of zero in $H^1_0$, with
\[
 \|f(u)\|_{H^1}\le C_f\|u\|_{H^1}^3,\qquad
 \|f'(u)v\|_{H^1}\le C_f\|u\|_{H^1}^2\|v\|_{H^1}.
\]
For a positive mode the coordinate gradients are
\[
 \partial_{q_n}P=\partial_{\widetilde q_n}P
       =\frac1{\sqrt{2\mu_n}}\int_0^\pi f(u)\phi_n\dd x.
\]
The spectral coefficients of an $H^1_0$ function belong to $h^1$.
Writing $b_n=\int_0^\pi f(u)\phi_n\dd x$, we have
\[
 \sum_{n\in\Ell}\langle n\rangle^3
              \left|\frac{b_n}{\sqrt{\mu_n}}\right|^2
 \le C_\omega\sum_{n\in\Ell}\langle n\rangle^2|b_n|^2
 \le C_\omega\|f(u)\|_{H^1}^2.
\]
Thus the coordinate gradients, and hence both infinite components of
$X_P$, satisfy the first bound in \eqref{eq:nonlinear-smooth}.  The
finite hyperbolic components are bounded linear functionals of $f(u)$.
Differentiating replaces $f(u)$ by $f'(u)\,Du(w)[v]$, which gives the
second bound.  The same bounded linear maps applied term by term to
the convergent power series of $f(u)$ prove the asserted analyticity
and convergence of the homogeneous vector fields.  Termwise
integration in $x$ gives the scalar expansion.
\end{proof}

Fix a finite $J\subset\Ell$ and let $\mathcal N_+=\Ell\setminus J$.
Tangential variables are $(q_j,\widetilde q_j)_{j\in J}$; all variables indexed
by $\mathcal N_+$ or $\Hyp$ are normal.  The normal degree of a monomial
counts all its normal factors, with multiplicities.  Write
$I_a=q_{j_a}\widetilde q_{j_a}$.
We denote the corresponding normal subspace of $\cE^p$ by $\cE_N^p$.

\subsection{Divisor and matrix estimates}
For fixed total degree, only finitely many tangential combinations
$C=\sum_{j\in J}c_j\mu_j$ occur in a homological equation of normal
degree at most two.  The elliptic divisors are
\begin{equation}\label{eq:divisor-types}
 C,\qquad C\pm\mu_n,\qquad C\pm(\mu_n+\mu_m),\qquad
 C+\mu_n-\mu_m.
\end{equation}
For a nonzero tangential coefficient vector, \eqref{eq:affine-nr}
gives $C\notin\Z$; the zero vector gives $C=0$.
The purely tangential divisor vanishes when $C=0$.  Among the
two-mode difference divisors, the resonant case is exactly $C=0$,
$n=m$.  The remaining divisors are nonzero by \eqref{eq:affine-nr}.
We first estimate them and then prove that the corresponding
generating Hamiltonians have one-smoothing vector fields.

\begin{lemma}[Estimates for the divisors]\label{lem:wave-divisors}
Fix a finite collection $\mathcal C$ of integer tangential combinations
$C=\sum_{j\in J}c_j\mu_j$, $c_j\in\Z$.  Each is either zero or
belongs to $\R\setminus\Z$ by \eqref{eq:affine-nr}.
There is $c>0$, depending on the path, $J$, and $\mathcal C$, such that
\begin{align}
 |C\pm\mu_n|&\ge c\langle n\rangle,\label{eq:wave-one-div}\\
 |C\pm(\mu_n+\mu_m)|&\ge c(1+n+m),\label{eq:wave-sum-div}\\
 |C+\mu_n-\mu_m|&\ge c(1+|n-m|),\label{eq:wave-difference-div}
\end{align}
for $n,m\in\mathcal N_+$ and $C\in\mathcal C$, excluding the resonant
case $C=0$, $n=m$ in the last line.
\end{lemma}
\begin{proof}
The first two bounds follow from $\mu_n=n+o(1)$ once the largest index
is sufficiently large.  The remaining cases are finite and have
nonzero divisors by \eqref{eq:affine-nr}.

For differences write $\mu_n=n+r_n$, with $r_n\to0$.
If $C\ne0$, then
\[
 d_C:=\inf_{k\in\Z}\frac{|k+C|}{1+|k|}>0:
\]
the quotient tends to one as $|k|\to\infty$, and none of its finitely
many remaining values vanishes.  When both indices are sufficiently
large, $|r_n|+|r_m|\le d_C/2$, and hence
\[
 |n-m+C+r_n-r_m|\ge\tfrac12d_C(1+|n-m|).
\]
For $C=0$, $n\ne m$, use instead
$|n-m|\ge(1+|n-m|)/2$ and choose the tail so that
$|r_n|+|r_m|\le1/2$.
If one index is outside the chosen tail, it belongs to a finite set.
For each such index, the quotient in \eqref{eq:wave-difference-div}
tends to one as the other index tends to infinity.  The remaining
finite block is controlled by exact nonresonance.  Taking the minimum
over the finite collection $\mathcal C$ completes the proof.
\end{proof}

To estimate the divided matrices, we use the specific form of the
divisors.  A bound on their reciprocals alone would not suffice, since
entrywise multiplication by a bounded array need not preserve
boundedness of an operator.

\begin{lemma}[Estimates for the divided matrices]\label{lem:Schur}
Suppose $M$ and its transpose define bounded maps $h^{1/2}(\mathcal N_+)$
to $h^{3/2}(\mathcal N_+)$.  For any of the finitely many constants $C$
above, entrywise division by the nonresonant two-mode divisors in
\eqref{eq:divisor-types} gives a matrix $K$ such that both $K$ and
$K^{\mathsf T}$ have the same mapping property.  For the difference divisor with $C=0$, the resonant diagonal
entries are set to zero.  Writing $\|M\|_{\rm s}$ for the sum of the operator norms of $M$
and $M^{\mathsf T}$ from $h^{1/2}$ to $h^{3/2}$, one has
\[
 \|K\|_{\rm s}\le C_{\omega,J,\mathcal C}\|M\|_{\rm s},
\]
where $\mathcal C$ is the finite collection of tangential constants.
\end{lemma}
\begin{proof}
\noindent\textit{Step 1: difference divisors on the tail.}
Let $\mathfrak B=\mathcal L(h^{1/2},h^{3/2})$, initially restricted
to a cofinite set of indices.  The norm of an operator in this space is
the $\ell^2$ operator norm of its weighted matrix
$\langle n\rangle^{3/2}M_{nm}\langle m\rangle^{-1/2}$.
All diagonal groups used below commute with these weights.
On a sufficiently high tail,
write $\mu_n=n+r_n$, where $r_n\to0$.
Let $U(t)=\diag(e^{int})$, an isometry on every $h^p$.
If $C\notin\Z$, the operator
\begin{equation}\label{eq:Fourier-divider}
 T_C(M)=\frac{i}{e^{2\pi iC}-1}\int_0^{2\pi}
                 e^{iCt}U(t)MU(-t)\dd t
\end{equation}
has entries $M_{nm}/(n-m+C)$ and norm at most
$2\pi|e^{2\pi iC}-1|^{-1}\|M\|_{\cL(h^{1/2},h^{3/2})}$.
For each $v\in h^{1/2}$, strong continuity of $U(t)$ on the source
and target spaces makes the integrand applied to $v$ a continuous
$h^{3/2}$-valued function.  Its Bochner integral satisfies
\[
 \|T_C(M)v\|_{3/2}
 \le \frac{2\pi}{|e^{2\pi iC}-1|}\|M\|_{\mathfrak B}\|v\|_{1/2}.
\]
The $(n,m)$ entry follows by applying this integral to $e_m$ and
integrating $e^{i(n-m+C)t}$.  Thus $T_C$ is a bounded linear operator on $\mathfrak B$.

For $C=0$, use
\begin{equation}\label{eq:sawtooth}
 T_0(M)=\frac{i}{2\pi}\int_0^{2\pi}(t-\pi)U(t)MU(-t)\dd t.
\end{equation}
For $k\in\Z\setminus\{0\}$,
\[
 \frac{i}{2\pi}\int_0^{2\pi}(t-\pi)e^{ikt}\dd t=\frac1k,
 \qquad \int_0^{2\pi}(t-\pi)\dd t=0.
\]
Thus its diagonal entries vanish and its off-diagonal entries are
$M_{nm}/(n-m)$.  The same estimate for the Bochner integral gives
$\|T_0\|_{\mathcal L(\mathfrak B)}\le\pi/2$.
When $C=0$, work in the closed subspace of $\mathfrak B$ with zero
diagonal, after deleting the resonant diagonal of $M$.
This deletion is bounded: application to $e_n$ gives
$\langle n\rangle|M_{nn}|\le\|M\|_{\mathfrak B}$, and therefore
$\|\diag M\|_{\mathfrak B}\le\|M\|_{\mathfrak B}$.

Put $R=\diag(r_n)$ on the tail and $\mathcal A(K)=RK-KR$.
Then $\|\mathcal A\|\le2\sup_{n\ge N}|r_n|$ as an operator on
$\cL(h^{1/2},h^{3/2})$.  Choose $N$ so large that
$2\|T_C\|\sup_{n\ge N}|r_n|<1$.  The equation
\[
 K=T_C(M)-T_C\mathcal A(K)
\]
has the solution
\[
 K=\sum_{j\ge0}(-T_C\mathcal A)^jT_C(M),\qquad
 \|K\|_{\mathfrak B}
 \le\frac{\|T_C\|}{1-2\|T_C\|\sup_{n\ge N}|r_n|}
                                    \|M\|_{\mathfrak B}.
\]
For $C=0$ the commutator $\mathcal A$ preserves the zero-diagonal
subspace.  Taking matrix entries
shows
\[
 K_{nm}=\frac{M_{nm}}{n-m+C+r_n-r_m}
       =\frac{M_{nm}}{\mu_n-\mu_m+C}
\]
for every nonresonant entry; when $C=0$, $K_{nn}=0$.
Transposition replaces $C$ by $-C$ and changes an overall sign.
Applying the same argument to $M^{\mathsf T}$ gives the transpose bound.

\medskip
\noindent\textit{Step 2: sum divisors on the tail.}
For sum divisors, choose $N$ with $2\mu_N+C>0$ and write
$L=\diag(\mu_n)$ on the tail.  Then
\begin{equation}\label{eq:Laplace-divider}
 K=\int_0^\infty e^{-Ct}e^{-tL}Me^{-tL}\dd t
\end{equation}
has entries $M_{nm}/(\mu_n+\mu_m+C)$ and norm at most
$(2\mu_N+C)^{-1}\|M\|$.  The two semigroups are diagonal contractions
up to the common factor $e^{-t\mu_N}$ on the respective weighted
spaces.  This also proves the transpose estimate.  The other overall
sign is absorbed by replacing $C$ with $-C$.

\medskip
\noindent\textit{Step 3: rows and columns outside the tail.}
It remains to include the finitely many rows or columns with an index
below $N$.  For each such row or column the other index tends to
infinity, and the modulus of its divisor tends to infinity.
Its reciprocal sequence is bounded: only finitely many remaining
entries require exact nonresonance.  A fixed column belongs to $h^{3/2}$, whereas a fixed row is a
bounded functional on $h^{1/2}$.  Multiplication of their coefficients
by a bounded sequence preserves the respective norms.  In particular,
if $d_{nm}$ is the relevant divisor and $m$ is fixed, then
\[
 \left\|\left(\frac{M_{nm}}{d_{nm}}\right)_n\right\|_{3/2}
 \le\sup_n|d_{nm}^{-1}|\,\|Me_m\|_{3/2}
 \le C_m\|M\|_{\mathfrak B}\langle m\rangle^{1/2}.
\]
A row is estimated in the weighted dual norm
$\bigl(\sum_m\langle m\rangle^{-1}|M_{nm}|^2\bigr)^{1/2}$.
The finitely many constants are absorbed in $C_{\omega,J,\mathcal C}$.
The sum of these finitely many row and column operators is bounded,
and the same argument applied to $M^{\mathsf T}$ gives the transpose
bound.  Together with the tail estimates, this proves the lemma.
\end{proof}

\paragraph{The quartic coefficients and the divided operators.}
The quartic term in \eqref{eq:nonlinear-Hamiltonian} is
$P_4=(\kappa/4)\int u^4$.  We express its homological equation in terms of the operators
in \cref{lem:Schur}.  Restricting first to positive modes,
\begin{equation}\label{eq:wave-quartic-expansion}
 P_4^{\rm ell}=\frac\kappa{16}\sum_{i,j,k,l\in\Ell}
 \frac{\int_0^\pi\phi_i\phi_j\phi_k\phi_l\dd x}
      {\sqrt{\mu_i\mu_j\mu_k\mu_l}}
 (q_i+\widetilde q_i)(q_j+\widetilde q_j)
 (q_k+\widetilde q_k)(q_l+\widetilde q_l).
\end{equation}
The expansion holds for finitely supported sequences and extends to
the energy space through the continuous quartic integral.  Under a sign
change $\phi_n\mapsto s_n\phi_n$, $s_n\in\{\pm1\}$, the integral
coefficient acquires the factor $s_i s_j s_k s_l$.  The simultaneous
coordinate changes $(q_n,\widetilde q_n)\mapsto
(s_nq_n,s_n\widetilde q_n)$ leave the displacement and the Hamiltonian
unchanged.

Fix three tangential indices $a,b,c\in J$ and a choice of $q$ or
$\widetilde q$ at each index.  A contribution with one
normal factor is, up to a fixed scalar, a tangential cubic monomial
times $\sum_n f_n z_n$ or $\sum_n f_n\widetilde z_n$, where
\[
 f_n=\frac1{\sqrt{\mu_n}}\int_0^\pi
                                  \phi_a\phi_b\phi_c\phi_n\dd x.
\]
The product $\phi_a\phi_b\phi_c$ belongs to $H^1_0$, so $f\in h^{3/2}$
by the estimate in \cref{lem:smoothing}.  The homological equation
replaces $f_n$ by $-f_n/[i(C\pm\mu_n)]$.  By
\cref{lem:wave-divisors}, the divided sequence still belongs to
$h^{3/2}$.  This sequence is the normal coordinate gradient of the
corresponding term in the generating Hamiltonian.

For two fixed tangential indices $a,b\in J$, set $g=\phi_a\phi_b$ and
\begin{equation}\label{eq:wave-multiplication-matrix}
 M^g_{nm}=\frac{\int_0^\pi g\phi_n\phi_m\dd x}{\sqrt{\mu_n\mu_m}},
 \qquad n,m\in\mathcal N_+.
\end{equation}
Let $\mathcal U_N z=\sum_{n\in\mathcal N_+}z_n\phi_n$ and let $L_N$
be the restriction of $L$ to these modes.  Then
\[
 M^g=L_N^{-1/2}\mathcal U_N^*\,\mathsf M_g\,
                         \mathcal U_N L_N^{-1/2},
 \qquad \mathsf M_g f=gf.
\]
Here $\mathcal U_N^*$ is the $L^2$ adjoint of $\mathcal U_N$.  Moreover,
$\mathcal U_NL_N^{-1/2}$ maps $h^{1/2}$ into $H^1_0$,
$\mathsf M_g$ is bounded on $H^1_0$, and
$L_N^{-1/2}\mathcal U_N^*$ maps $H^1_0$ into $h^{3/2}$.  Consequently
\[
 \|M^g z\|_{3/2}\le C_\omega\|g\|_{H^1}\|z\|_{1/2}.
\]
The matrix is symmetric, so its transpose has the same bound.
Depending on the two normal coordinate types, the homological
equation replaces $M^g_{nm}$ by
$-M^g_{nm}/[i(C\pm(\mu_n+\mu_m))]$ or
$-M^g_{nm}/[i(C+\mu_n-\mu_m)]$, up to the fixed scalar factors
from the tangential coordinates.  For the difference divisor, the
resonant diagonal with $C=0$ is omitted.
These are the divided matrices of \cref{lem:Schur}.
The corresponding quadratic forms have gradients given by these
matrices and their transposes applied to the normal variables.

At higher degrees the matrices produced by Poisson brackets need not
have the multiplication form \eqref{eq:wave-multiplication-matrix}.
The general operator estimate in \cref{lem:Schur} therefore applies
throughout the finite normalization.  We now solve the homological
equation in this class of Hamiltonians.

\begin{lemma}[The homological equation]\label{lem:homological}
Let $P$ be a real homogeneous polynomial Hamiltonian of even degree
$d\ge4$ whose vector field maps $\cE^{1/2}$ to $\cE^{3/2}$ continuously.
Let $P^{\le2}$ be its part of normal degree at most two.
Its resonant part $Z$ consists precisely of terms of the form
\begin{equation}\label{eq:resonant-structure}
 Z_T(I)+\sum_{n\in\mathcal N_+}Z_n(I)q_n\widetilde q_n
                  +\sum_{h\in\Hyp}Y_h(I)a_hb_h,
\end{equation}
where $Z_T$ is homogeneous of action degree $d/2$ and $Z_n,Y_h$ have
action degree $d/2-1$.
There is a real homogeneous generating Hamiltonian $F$ of degree $d$,
containing only terms of normal degree at most two, whose vector field
has the same smoothing property and satisfies
\begin{equation}\label{eq:homological}
 \{H_2,F\}+P^{\le2}-Z=0.
\end{equation}
The identity holds as an identity of continuous scalar Hamiltonians on
$\cE^{1/2}$.  If
\[
 \|P\|_d:=\sup_{\|w\|_{1/2}\le1}\|X_P(w)\|_{3/2},
\]
then, for a constant depending only on the fixed path, $J$, and $d$,
\begin{equation}\label{eq:homological-bound}
 \|F\|_d+\|Z\|_d+\|P^{\le2}\|_d\le C_{\omega,J,d}\|P\|_d.
\end{equation}
\end{lemma}
\begin{proof}
\noindent\textit{Step 1: estimates for the coefficient vectors and matrices.}
Write $w=(t,\zeta)$, where $t$ denotes the tangential coordinates and
$\zeta$ the normal coordinates.  Expanding $P(t,s\zeta)$ in the scalar
$s$ defines the homogeneous components $P^{(r)}$ of normal degree $r$.
The tangential gradient of $P^{(r)}$ has normal degree $r$, whereas its
normal gradient has normal degree $r-1$.  They are the coefficients of $s^r$ and $s^{r-1}$, respectively,
in the corresponding components of $X_P(t,s\zeta)$.
Cauchy estimates in $s$ bound these coefficients in the respective
vector-field norms.  Polarization gives
the equivalent bounded multilinear estimates, with constants depending
only on $d$.  Thus the projection to normal degree at most two preserves
smoothing, including both parts of the gradient.

Next expand in the finitely many tangential monomials.  The terms of
normal degrees one and two have the forms
\[
 T_{d-1}(t)\ell_f(\zeta_1),\qquad
 T_{d-2}(t)\mathcal B_M(\zeta_1,\zeta_2),
\]
where $T_j$ is a tangential monomial of degree $j$,
$\ell_f(z)=\sum_n f_nz_n$, and
$\mathcal B_M(z,v)=\sum_n(Mz)_n v_n$.
For a quadratic term in two identical coordinate types, choose
$M=M^{\mathsf T}$; only this symmetric part contributes to the
Hamiltonian.  Its derivative is $2Mz$.
For two different coordinate types, the gradients with respect to
$v$ and $z$ are $Mz$ and $M^{\mathsf T}v$, respectively.
Coefficient extraction in the finitely many tangential variables
therefore gives $f\in h^{3/2}$ and bounded maps
$M,M^{\mathsf T}:h^{1/2}\to h^{3/2}$,
with norms at most $C_{J,d}\|P\|_d$.  These are the estimates required
by \cref{lem:Schur}.  Terms involving hyperbolic coordinates have
a finite-dimensional domain or range and satisfy the same coefficient
estimates.

\medskip
\noindent\textit{Step 2: identification of the resonant monomials.}
For a monomial with elliptic exponents $\alpha,\beta$ and hyperbolic
exponents $r,s$, the bracket with $H_2$ multiplies it by
\begin{equation}\label{eq:hom-eigenvalue}
 i\sum_{n\in\Ell}(\alpha_n-\beta_n)\mu_n
                    +\sum_{h\in\Hyp}(s_h-r_h)\chi_h.
\end{equation}
If no hyperbolic variable occurs, \eqref{eq:affine-nr} with $a=0$
implies that a resonant monomial has equal exponents of $q_n$ and
$\widetilde q_n$ at each mode.  This gives the first two types in
\eqref{eq:resonant-structure}.  If one hyperbolic variable occurs, the
real part of \eqref{eq:hom-eigenvalue} is nonzero.  With two such variables,
their real contributions can cancel only for $a_hb_h$ at the same index:
the positive numbers $\chi_h$ are distinct because the eigenvalues are
simple.  The tangential contribution must then also vanish, leaving
exactly the last type in \eqref{eq:resonant-structure}.

\medskip
\noindent\textit{Step 3: division and the vector-field estimates.}
We divide each nonresonant coefficient by the negative of
\eqref{eq:hom-eigenvalue}.  With no normal variable this is a finite
calculation.  For one elliptic normal variable, division multiplies the
coefficient vector $f\in h^{3/2}$ by the bounded sequence $(C\pm\mu_n)^{-1}$.
For two elliptic normal variables, \cref{lem:Schur} bounds the divided
matrix and its transpose from $h^{1/2}$ to $h^{3/2}$.

One hyperbolic and one elliptic variable give a vector coefficient and
a divisor with nonzero real part $\pm\chi_h$.  Its reciprocal is
bounded uniformly in the elliptic index.  Terms whose normal variables
are all hyperbolic involve only finitely many coefficients, and each
nonresonant divisor is nonzero by Step~2.  To obtain the vector-field
bound from these coefficient estimates, put
$A=|t|$ and $R=\|\zeta\|_{1/2}$, with both elliptic coordinate types
included.  For a divided coefficient $f\in\cE_N^{3/2}$,
\[
 \begin{aligned}
 |\partial_t(T_{d-1}\ell_f)|
       &\le C A^{d-2}\|f\|_{3/2}R,\\
 \|\partial_\zeta(T_{d-1}\ell_f)\|_{3/2}
       &\le C A^{d-1}\|f\|_{3/2}.
 \end{aligned}
\]
For a divided matrix $K$, the analogous estimates are
\[
 \begin{aligned}
 |\partial_t(T_{d-2}\mathcal B_K)|
       &\le C A^{d-3}\|K\|_{\rm s}R^2,\\
 \|\partial_\zeta(T_{d-2}\mathcal B_K)\|_{3/2}
       &\le C A^{d-2}\|K\|_{\rm s}R.
 \end{aligned}
\]
They follow from the product rule and
$|\mathcal B_K(z,v)|\le C\|K\|_{\rm s}\|z\|_{1/2}\|v\|_{1/2}$.
The coordinate gradients give the vector-field components by
multiplication by $\pm i$, or by interchange and a sign in the finite
hyperbolic pairs.  Tangential weighted norms are equivalent to
Euclidean norms for fixed $J$.  Since $A+R\le C_J\|w\|_{1/2}$,
summing these estimates bounds the entire vector field by
$C\|P\|_d\|w\|_{1/2}^{d-1}$.
The divided-matrix bounds in \cref{lem:Schur} and the finite scalar
divisors give the $F$ estimate in \eqref{eq:homological-bound}.
The resonant diagonal projection is bounded as well: if
$M:h^{1/2}\to h^{3/2}$ is bounded, then
$|M_{nn}|\le C\langle n\rangle^{-1}$, so its diagonal defines a map
of the same class.  Tangential coefficient projections are finite, and the finite
hyperbolic projections are bounded.  This proves the $Z$ estimate;
the estimate for $P^{\le2}$ follows from Step~1.

\medskip
\noindent\textit{Step 4: reality and the homological identity.}
The constructed series define continuous polynomial Hamiltonians:
linear normal terms use bounded functionals, quadratic terms use the
bounded bilinear forms above.  Their finite-coordinate truncations converge at every point of the
energy space, and the extension is unique by density.  To verify the
reality conditions, consider the antilinear involution
\[
 \mathcal R(q,\widetilde q,a,b)
       =(\overline{\widetilde q},\overline q,\overline a,\overline b).
\]
Its fixed points correspond to real-valued wave data.  Define
$G^\#(w)=\overline{G(\mathcal Rw)}$.  The reality of $H_2$
and $P$ means $H_2^\#=H_2$ and $P^\#=P$, and the Poisson bracket
satisfies $\{G,F\}^\#=\{G^\#,F^\#\}$.  Hence $F^\#$ solves the
same homological equation as $F$ and has zero resonant part.
Uniqueness of the divided coefficients gives $F^\#=F$, or
$F(\mathcal Rw)=\overline{F(w)}$.
Consequently $X_F$ preserves the reality conditions.

Finally \eqref{eq:homological} holds first on finite sequences by
\eqref{eq:hom-eigenvalue}.  The differential of $H_2$ is continuous on
the energy space, with
$|\dd H_2(w)[v]|\le C\|w\|_{1/2}\|v\|_{1/2}$.
Thus $\dd H_2[X_F]$ is continuous there.  Density extends the identity
to the whole energy space.
\end{proof}

\subsection{Construction of the normal form}
We use successive Lie transforms, given by the time-one maps of
homogeneous generating Hamiltonians.  For the general normal-form
construction for Hamiltonian PDEs, see
\cite{Bambusi2003,BambusiGrebert2006}.

\begin{proposition}[Partial Birkhoff normal form up to degree ten]\label{prop:NF}
There is a real-analytic canonical transformation $\Phi$ near zero in
$\cE^{1/2}$, with
\begin{equation}\label{eq:Phi-bounds}
 \|\Phi(w)-w\|_{3/2}\le C\|w\|_{1/2}^3,
 \qquad \|D\Phi(w)-\Id\|_{\cL(\cE^{1/2},\cE^{3/2})}
                                      \le C\|w\|_{1/2}^2,
\end{equation}
such that
\begin{equation}\label{eq:NF}
 H\circ\Phi=H_2+Z_T(I)+\sum_{n\in\mathcal N_+}Z_n(I)q_n\widetilde q_n
       +\sum_{h\in\Hyp}Y_h(I)a_hb_h+P_{\ge3}+R_{12}.
\end{equation}
Here $Z_T$ has action degrees from two to five, while $Z_n$ and $Y_h$
have action degrees from one to four, and $P_{\ge3}$ is a sum of homogeneous
polynomials of degrees $4,6,8,10$, each of normal degree at least three.
All the nonlinear terms have analytic vector fields from
$\cE^{1/2}$ to $\cE^{3/2}$.  On a fixed smaller ball,
\begin{equation}\label{eq:R12}
 \|X_{R_{12}}(w)\|_{3/2}\le C\|w\|_{1/2}^{11},\qquad
 \|DX_{R_{12}}(w)\|_{\cL(\cE^{1/2},\cE^{3/2})}
                                      \le C\|w\|_{1/2}^{10}.
\end{equation}
The coefficients of $Z_n$ are $O_{\omega,f,J}(\langle n\rangle^{-1})$,
uniformly in the finitely many action monomials.
\end{proposition}
\begin{proof}
We use the Hamiltonian degree throughout the induction.  The energy
space is $\cE^{1/2}$; every nonlinear vector field is estimated with
values in $\cE^{3/2}$.
By \cref{lem:smoothing}, the initial nonlinear Hamiltonian is even and
analytic, with Taylor decomposition
\[
 P=P_4+P_6+P_8+P_{10}+R_{12}^{\rm init}.
\]
The original homogeneous terms are those in
\eqref{eq:nonlinear-Hamiltonian}.  On a sufficiently small ball, the Taylor remainder satisfies
\begin{equation}\label{eq:initial-remainder}
 \|X_{R_{12}^{\rm init}}(w)\|_{3/2}\le C\|w\|_{1/2}^{11},\qquad
 \|DX_{R_{12}^{\rm init}}(w)\|_{\cL(\cE^{1/2},\cE^{3/2})}
                                      \le C\|w\|_{1/2}^{10}.
\end{equation}
These bounds follow by summing the convergent homogeneous series after
its terms of vector-field degrees $3,5,7,9$ have been removed.

\medskip
\noindent\textit{Step 1: flows of homogeneous generating Hamiltonians.}
Let $F_d$ be a homogeneous Hamiltonian of even degree $d\ge4$ with
$\|X_{F_d}(w)\|_{3/2}\le C\|w\|_{1/2}^{d-1}$.  Differentiating its
bounded multilinear representation gives
\[
 \|DX_{F_d}(w)\|_{\cL(\cE^{1/2},\cE^{3/2})}
       \le C_d\|w\|_{1/2}^{d-2}.
\]
The local existence theorem for analytic vector fields and the inclusion
$\cE^{3/2}\hookrightarrow\cE^{1/2}$ give a local flow
$\varphi_d^t$.  The integral equation and a continuation argument on a sufficiently
small ball give
\begin{equation}\label{eq:flow-estimates}
 \begin{split}
 \|\varphi_d^t(w)\|_{1/2}&\le2\|w\|_{1/2},\\
 \|\varphi_d^t(w)-w\|_{3/2}&\le C\|w\|_{1/2}^{d-1},\qquad |t|\le1.
 \end{split}
\end{equation}
The second inequality follows by integrating $X_{F_d}$ in
$\cE^{3/2}$.  In the variational equation
\[
 D\varphi_d^t(w)=\Id+
   \int_0^t DX_{F_d}(\varphi_d^s(w))D\varphi_d^s(w)\dd s,
\]
Gronwall's inequality first bounds the differential on $\cE^{1/2}$.
Using the bound for $DX_{F_d}$ with values in $\cE^{3/2}$ in the
same equation gives
\begin{equation}\label{eq:flow-differential}
 \|D\varphi_d^t(w)-\Id\|_{\cL(\cE^{1/2},\cE^{3/2})}
       \le C\|w\|_{1/2}^{d-2}.
\end{equation}
On a smaller ball, the differential and its inverse are uniformly
bounded on both spaces.  On $\cE^{3/2}$ this follows from the Neumann series, since
$D\varphi_d^t-\Id$ has small operator norm on that space.
The inverse flow is $\varphi_d^{-t}$ on its image.

The weak symplectic form \eqref{eq:symp} is continuous and
nondegenerate.  For fixed tangent vectors $U,V$, the variational
equation and $\Om(X_{F_d},\cdot)=\dd F_d$ give
\[
 \frac{\dd}{\dd t}
 \Om(D\varphi_d^tU,D\varphi_d^tV)
 =D^2F_d[D\varphi_d^tU,D\varphi_d^tV]
  -D^2F_d[D\varphi_d^tV,D\varphi_d^tU]=0.
\]
Thus $\varphi_d^t$ is symplectic.  Uniqueness of the ODE preserves the
reality conditions, and the flow has a holomorphic extension
to the complex energy ball.

\medskip
\noindent\textit{Step 2: brackets and pullbacks preserve smoothing.}
For nonlinear Hamiltonians $F,G$ with analytic vector fields from
$\cE^{1/2}$ to $\cE^{3/2}$, the sign
convention in \eqref{eq:symp} gives
\begin{equation}\label{eq:bracket-field}
 X_{\{G,F\}}=DX_G X_F-DX_F X_G.
\end{equation}
This identity follows by pairing with an arbitrary tangent vector and
using symmetry of the Hessians.  Both compositions take values in
$\cE^{3/2}$, since their arguments in that space also belong to
$\cE^{1/2}$.  If $G$ and $F$ are homogeneous of degrees $e$ and $d$,
respectively, polarization bounds the right-hand side by
$C_{d,e}\|G\|_e\|F\|_d\|w\|_{1/2}^{e+d-3}$.
Hence the scalar bracket has degree $e+d-2$ and has a smoothing vector
field.  For analytic Hamiltonians, the corresponding estimates follow
on smaller balls by Cauchy estimates.

Symplecticity also gives the exact pullback formula
\begin{equation}\label{eq:pullback-field}
 X_{G\circ\varphi_d^t}(w)
   =(D\varphi_d^t(w))^{-1}X_G(\varphi_d^t(w)).
\end{equation}
Indeed, pairing its right-hand side with $V$ gives
$\dd G(\varphi_d^t(w))[D\varphi_d^t(w)V]$.
The bound for $(D\varphi_d^t)^{-1}$ on $\cE^{3/2}$ from Step~1
shows that this vector
field is analytic into $\cE^{3/2}$ on a common smaller ball, uniformly
for $|t|\le1$.

\medskip
\noindent\textit{Step 3: one normalization step.}
Suppose the current Hamiltonian is $H_2+P$, and let $P_d$ denote the
full homogeneous component of degree $d$ in $P$.  Put
\[
 Q_d=P_d^{\le2}-Z_d,
 \qquad \{H_2,F_d\}=-Q_d,
\]
using \cref{lem:homological}.  The scalar functional $H_2$ is continuous
quadratic on the energy space.  The chain rule along the bounded flow
$\varphi_d^t$ is therefore valid and gives
$\frac{\dd}{\dd t}(H_2\circ\varphi_d^t)=-Q_d\circ\varphi_d^t$.
Integrating this identity twice, and the nonlinear one once, yields
\begin{align}
 H_2\circ\varphi_d^1
   &=H_2-Q_d-
       \int_0^1(1-t)\{Q_d,F_d\}\circ\varphi_d^t\dd t,
       \label{eq:linear-pullback}\\
 P\circ\varphi_d^1
   &=P+\int_0^1\{P,F_d\}\circ\varphi_d^t\dd t.
       \label{eq:nonlinear-pullback}
\end{align}
All brackets inside the integrals involve bounded smoothing vector
fields.  By Step~2 their pullbacks, and the integrals of those
pullbacks, have analytic smoothing vector fields.  Thus both pullback formulas are valid on the energy space.

The integral in \eqref{eq:linear-pullback} has no Taylor terms of
degree below $2d-2$.  Since $P$ has no terms of degree below four,
the integral in \eqref{eq:nonlinear-pullback} has no Taylor terms of
degree below $d+2$.
Consequently the degree-$d$ term after transformation is exactly
\[
 P_d-Q_d=Z_d+P_d^{\ge3},
\]
and all lower total degrees are unchanged.  At the first step, for
example, the degree-six term becomes
\begin{equation}\label{eq:degree-six-recursion}
 P_6^{\rm new}=P_6+\{P_4,F_4\}-\tfrac12\{Q_4,F_4\}.
\end{equation}
Here $P_4$ and $P_6$ are the original homogeneous terms, including all
normal degrees; in particular $P_6=(a_5/6)\int u^6$ participates in the
next homological equation.  More generally, write
$\operatorname{ad}_F G=\{G,F\}$.
For each homogeneous component of degree $e\le10$, the finite Taylor
expansion along the flow gives
\[
 [H\circ\varphi_d^1]_e
  =\sum_{j\ge0}\frac1{j!}
     \operatorname{ad}_{F_d}^{j}
                   [H]_{e-j(d-2)},\qquad e\le10,
\]
where absent or negative degrees contribute zero and the sum is finite.
Brackets with the quadratic term are evaluated using
$\{H_2,F_d\}=-Q_d$; all subsequent brackets have bounded smoothing
vector fields.  The integral formulas above justify the remainder in this finite
Taylor expansion.

\medskip
\noindent\textit{Step 4: finite induction and the remainder.}
Starting with $H_2+P$, perform Step~3 successively for $d=4,6,8,10$.
At each step, $P_d$ means the homogeneous term of the Hamiltonian
already transformed by all preceding flows.  Thus brackets involving
previously retained terms are included when computing the next
homogeneous term.  Terms of normal degree at least three can contribute to subsequent
homological equations, since a Poisson bracket differentiates the
normal variables and may produce terms of normal degree two.
For instance, differentiating a normal variable in each of two terms
of normal degrees three and one gives normal degree $3+1-2=2$
in their Poisson bracket.  Thus the newly formed
degree-six term in \eqref{eq:degree-six-recursion} must be decomposed
again before solving for $F_6$.

The induction may be written as
\[
 H^{[d-2]}=H_2+
   \sum_{\substack{4\le e<d\\e\ \mathrm{even}}}(Z_e+W_e)
                   +P_d+O_{\rm Ham}(d+2),
\]
where each $W_e$ has normal degree at least three, and
$O_{\rm Ham}(d+2)$ has a smoothing vector field vanishing to order
$d+1$.  Pullback by $\varphi_d^1$ replaces $P_d$ by
$Z_d+W_d$, with $W_d=P_d^{\ge3}$, and changes only terms of higher
total degree.  This verifies the induction hypothesis at $d+2$.

Set
\[
 \Phi=\varphi_4^1\circ\varphi_6^1\circ
                       \varphi_8^1\circ\varphi_{10}^1.
\]
The finite intersection of suitably small energy balls gives a common
domain for this composition.  Its bounds \eqref{eq:Phi-bounds} follow
from \eqref{eq:flow-estimates}--\eqref{eq:flow-differential} and the
chain rule; the degree-four generating Hamiltonian determines the leading terms.
The original Hamiltonian is even because $f$ is odd.  Every generating
Hamiltonian has even degree and its vector field is odd.  Its flow is
therefore odd, so the transformed Hamiltonian remains
even.  After the degree-ten step, the remaining scalar Taylor terms have
degree at least twelve.

By \eqref{eq:linear-pullback}--\eqref{eq:nonlinear-pullback}, the
nonlinear part of the final Hamiltonian has an analytic vector field
into $\cE^{3/2}$.  Its Taylor coefficients of degree at most nine are
the vector fields of the corresponding scalar homogeneous terms; this follows in each coordinate by differentiating
the Taylor expansion.  Subtracting them gives an analytic vector field that vanishes to
order eleven at the origin.  Taylor estimates on a fixed smaller ball
give the first bound in \eqref{eq:R12}.  The derivative bound follows
from its convergent multilinear expansion, or from a Cauchy estimate
on a ball of radius proportional to $\|w\|_{1/2}$; at the origin it
follows from the vanishing Taylor coefficients.

Finally, \cref{lem:homological} identifies every resonant polynomial
as \eqref{eq:resonant-structure}.  Collecting degrees $4,6,8,10$ gives
the stated action degrees and the polynomial $P_{\ge3}$.
If $M:h^{1/2}\to h^{3/2}$ is any coefficient matrix of a quadratic
normal term, application to the coordinate vector $e_n$ gives
\[
 \langle n\rangle^{3/2}|M_{nn}|
   \le\|Me_n\|_{3/2}\le\|M\|\langle n\rangle^{1/2}.
\]
Thus $|M_{nn}|\le\|M\|\langle n\rangle^{-1}$.
There are finitely many tangential monomials of degree at most ten,
which proves the uniform coefficient bound for $Z_n$.
\end{proof}

\subsection{Quartic normal-form coefficients}\label{sec:coefficients}
For the full set of positive modes, write
$\mathcal I_n=q_n\widetilde q_n$, so that $I_a=\mathcal I_{j_a}$.
The quartic action projection is
\begin{equation}\label{eq:quartic-actions}
 Z_4^{\rm ell}=
  \frac{3\kappa}{8}\sum_{n\in\Ell}\frac{C_{nn}}{\mu_n^2}\mathcal I_n^2
  +\frac{3\kappa}{2}\sum_{\substack{m<n\\m,n\in\Ell}}
                              \frac{C_{mn}}{\mu_m\mu_n}\mathcal I_m\mathcal I_n.
\end{equation}
To verify the coefficients, insert
$u_n=(q_n+\widetilde q_n)/\sqrt{2\mu_n}$ into $(\kappa/4)\int u^4$.
The coefficient of $q_n^2\widetilde q_n^2$ in $(q_n+\widetilde q_n)^4$ is $6$.
For two distinct modes there are $6$ placements of the two copies of
each mode and a factor $2\cdot2$ from their two quadratic expansions.
These multiplicities give \eqref{eq:quartic-actions}.
It follows that
\begin{equation}\label{eq:leading-coefficients}
 \begin{split}
 Z_T(I)&=\tfrac12 I^{\mathsf T}A_JI+O(|I|^3),\\
 Z_n(I)&=B_n\cdot I+O(\langle n\rangle^{-1}|I|^2),\qquad
 (B_n)_a=\frac{3\kappa C_{n,j_a}}{2\mu_n\mu_{j_a}},\\
 Y_h(I)&=D_h\cdot I+O(|I|^2),\qquad
 (D_h)_a=-\frac{3\kappa C_{h,j_a}}{2\chi_h\mu_{j_a}}.
 \end{split}
\end{equation}
The matrix $A_J$ is \eqref{eq:twist}.  For the last coefficient use
$(a_h-b_h)^2$ in the hyperbolic displacement: its $a_hb_h$ coefficient
is $-2$, while the elliptic quadratic factor contributes $2I_a$.
All remainders in \eqref{eq:leading-coefficients} are finite polynomials;
their derivatives on a fixed small action box satisfy the corresponding
uniform bounds.  In particular $\sup_n\langle n\rangle|B_n|<\infty$.

\section{KAM construction}\label{sec:KAM}
\subsection{An abstract KAM theorem}\label{sec:abstract-KAM}
We use the analytic KAM theorem of P\"oschel
\cite[Theorems A and B, Corollary C]{Poschel1996}, with exponential
sequence weight zero, $p=1/2$, $\bar p=3/2$, and spectral exponent $d=1$.
The case with finitely many hyperbolic pairs uses the analytic extension
in \cite[Theorem 2.5]{YuanHyperbolic2003}.  We state the hypotheses
in the notation of the normal form; the theorem and the associated
measure estimates are given in Appendix~\ref{app:KAM}.

Let $\zeta=(z,\widetilde z,a,b)$ denote all normal variables and let
$\|\zeta\|_p$ include the finite hyperbolic components.  Consider
\begin{equation}\label{eq:abstract-N}
 N_\xi=\omega(\xi)\cdot y+
       \sum_{n\in\mathcal N_+}\Omega_n(\xi)z_n\widetilde z_n
                       +\sum_{h\in\Hyp}\Lambda_h(\xi)a_hb_h
\end{equation}
on
\[
 D(s,r)=\{|\Im\theta|<s,\ |y|<r^2,\ \|\zeta\|_{1/2}<r\}.
\]
For a vector field define
\begin{equation}\label{eq:scaled-norm}
 \|X\|_r=\sup_{D(s,r)}
 \left(|X^\theta|+r^{-2}|X^y|+r^{-1}\|X^\zeta\|_{3/2}\right).
\end{equation}
The superscript $\Lip$ denotes the Lipschitz seminorm in $\xi$, with
phase variables fixed.  All finite-dimensional norms are equivalent;
we use the maximum norm on action boxes.

We verify the following frequency and regularity hypotheses:
\begin{enumerate}[label=\textup{(K\arabic*)}]
\item On a fixed compact reference box, $\omega$ is bi-Lipschitz onto
its image.  For $k\in\Z^b$ and finitely supported
$l\in\Z^{(\mathcal N_+)}$, put $|l|=\sum_n|l_n|$.
For every $(k,l)\ne0$ with $|l|\le2$, the resonance function
$k\cdot\omega+l\cdot\Omega$ has a Lebesgue-null zero set.  Pure normal
combinations with $1\le|l|\le2$ are nonzero throughout that box.
\item After increasing enumeration of the elliptic normal indices,
$\Omega_{n_r}=r+c+O(r^{-\eta})$ for some fixed $c$ and $0<\eta<1$,
with a Lipschitz remainder in the norm $\sup_r r^\eta|\cdot|$.
Furthermore
\begin{equation}\label{eq:KAM-gap}
 \frac{\Omega_{n_r}-\Omega_{n_s}}{r-s}=1+O(s^{-\eta})\qquad(r>s)
\end{equation}
on the tail, uniformly in the parameters.
\item The perturbation vector field is analytic from the phase space
with index $1/2$ to that with index $3/2$, and Lipschitz in the parameters.
If hyperbolic pairs are present, their rates $\Lambda_h$ are
Lipschitz in the parameters, positive, and pairwise separated
uniformly on the parameter box.
\end{enumerate}
These hypotheses give the specialization of the cited KAM theorems
used here.  The finitely many exceptional normal indices are included
in the frequency bounds, and the reindexed sequence norms are equivalent
to the original ones.  The decay exponent in their notation is
$\delta=-\eta$, and $p-\bar p=-1\le\delta<0$.

With uniform bounds $M\ge1$ on the frequency norms and $L\ge1$ on the
inverse Lipschitz constant, the analytic smallness condition is
\begin{equation}\label{eq:smallness}
 \epsilon_P:=\|X_P\|_r+\frac\alpha M\|X_P\|_r^{\Lip}
                  \le\gamma\alpha,\qquad 0<\alpha\le1.
\end{equation}
It yields analytic invariant tori and embedding errors $O(\epsilon_P/\alpha)$
in scaled phase norm.  We take the Diophantine exponent sufficiently
large, for example
\begin{equation}\label{eq:tau}
 \tau> (b+3)(1+\eta)/\eta,
\end{equation}
which satisfies the $d=1$ summability condition in
\cite[equation (22)]{Poschel1996}.
P\"oschel's Theorem~B singles out finitely many exceptional
index pairs $(k,l)$.  Let
$\mathcal R_{\mathrm{ne}}(\alpha)$ denote the union of the excluded zones
with nonexceptional index pairs over all KAM steps.  The measure estimate is
\begin{equation}\label{eq:geometric}
 |\mathcal R_{\mathrm{ne}}(\alpha)|
 \le C(\operatorname{diam}\Pi)^{b-1}\alpha^\vartheta,
 \qquad \vartheta=\frac\eta{1+\eta}.
\end{equation}
After reducing $\gamma$ as in Corollary~C, the exceptional index pairs occur
only in the initial exclusion.  In \cref{lem:uniform} we show that their
zones are empty on sufficiently small action boxes and obtain the same
bound for the entire excluded set.

The finite hyperbolic modes satisfy the hypotheses of Yuan's analytic
extension.  Their positive, pairwise separated rates imply that the
additional resonance zones are empty for sufficiently small $\alpha$;
see Appendix~\ref{app:hyperbolic}.  Thus \eqref{eq:geometric} remains
the measure estimate for the elliptic--hyperbolic problem.

\subsection{Verification of the frequency hypotheses}
Introduce tangential action--angle coordinates
\begin{equation}\label{eq:AA}
 q_{j_a}=\sqrt{\xi_a+y_a}\,e^{-i\theta_a},\qquad
 \widetilde q_{j_a}=\sqrt{\xi_a+y_a}\,e^{i\theta_a},
\end{equation}
and retain $(z_n,\widetilde z_n)=(q_n,\widetilde q_n)$ for $n\in\mathcal N_+$ and
$(a_h,b_h)$ for $h\in\Hyp$.  A direct calculation gives
$i\dd q_{j_a}\wedge\dd\widetilde q_{j_a}=\dd\theta_a\wedge\dd y_a$.
The frequencies from \eqref{eq:NF} are
\begin{equation}\label{eq:frequencies}
 \omega(\xi)=\mu_J+\nabla Z_T(\xi),\qquad
 \Omega_n(\xi)=\mu_n+Z_n(\xi),\qquad
 \Lambda_h(\xi)=\chi_h+Y_h(\xi).
\end{equation}
All higher action terms from the finite normal form are included here.

\begin{proposition}\label{prop:frequency-check}
There is a fixed box $\Pi^0=[0,2\nu_*]^b$ on which the frequencies
\eqref{eq:frequencies} satisfy \textup{(K1)} and \textup{(K2)}, and the
hyperbolic condition in \textup{(K3)}.  Their required Lipschitz bounds
and the inverse bound are finite and independent of the shrinking
scale $\nu$.
\end{proposition}
\begin{proof}
By \eqref{eq:leading-coefficients},
$D\omega(\xi)=A_J+O(|\xi|)$.  Write
$\omega(\xi)=\mu_J+A_J\xi+R(\xi)$ and choose $\nu_*$ so small that
$\Lip(R)\le(2\|A_J^{-1}\|)^{-1}$.  Then for $\xi,\xi'\in\Pi^0$,
\[
 |\omega(\xi)-\omega(\xi')|
 \ge\frac{|\xi-\xi'|}{2\|A_J^{-1}\|}.
\]
The upper Lipschitz bound follows from bounded derivatives.  This
proves bi-Lipschitz invertibility with uniform constants.

For fixed $(k,l)\ne0$, the elliptic resonance function is a polynomial
in $\xi$.  At $\xi=0$ its value is
\[
 k\cdot\mu_J+\sum_{n\in\mathcal N_+}l_n\mu_n\ne0
\]
by \eqref{eq:affine-nr}; the two index sets are disjoint.
It is therefore a nonzero polynomial, whose zero set has Lebesgue
measure zero.  This verifies the zero-set condition in \textup{(K1)}.

For pure normal combinations, all positive $\mu_n$ have a positive
minimum.  Their sums are bounded below by twice that minimum.
Moreover, consecutive positive gaps tend to $1$, so the finitely
many initial gaps and the tail give
\[
 \inf_{\substack{n,m\in\Ell\\n\ne m}}|\mu_n-\mu_m|>0.
\]
Thus $|l\cdot\mu|$ is uniformly positive for $1\le|l|\le2$.
The bound $|Z_n(\xi)|\le C\langle n\rangle^{-1}|\xi|$
keeps these combinations nonzero when $\nu_*$ is small.
The finitely many $\chi_h>0$ are distinct, so the same reduction of
$\nu_*$ keeps all $\Lambda_h$ positive and distinct.

Enumerate $\mathcal N_+=\{n_1<n_2<\cdots\}$.
Eventually $n_r=r+h(\omega)+b$.  Set $\bar\Omega_r=n_r$; then
$\bar\Omega_r=r+h(\omega)+b$ for all sufficiently large $r$.
\Cref{prop:spectral} gives $\mu_{n_r}-\bar\Omega_r=O(r^{-\eta})$.
The coefficient bounds in \cref{prop:NF} give
\begin{equation}\label{eq:normal-coeff-bounds}
 \sup_{\xi\in\Pi^0}\langle n\rangle
       \bigl(|Z_n(\xi)|+|DZ_n(\xi)|+|D^2Z_n(\xi)|\bigr)\le C.
\end{equation}
Consequently $\Omega_{n_r}-\bar\Omega_r$ is Lipschitz in the norm
$\sup_r r^\eta|\cdot|$.  For $r>s$ on the tail, the error in
$\Omega_{n_r}-\Omega_{n_s}-(r-s)$ is bounded by $Cs^{-\eta}$.
Division by $r-s\ge1$ gives \eqref{eq:KAM-gap}.
The finitely many remaining indices are included by increasing the constants.  This proves the assertion.
\end{proof}

Choose $M,L$ to bound the frequency norms and the inverse Lipschitz
constant on the reference box, including the
spectral remainder $\mu_n-n$ and the finite hyperbolic rates.  These bounds hold
on $\Pi^0$ and hence on every $\Pi_\nu=[\nu,2\nu]^b\subset\Pi^0$.

\begin{lemma}[Uniformity on shrinking boxes]\label{lem:uniform}
Let $\alpha_\nu>0$ with $\alpha_\nu\to0$ as $\nu\downarrow0$.
There are constants $\gamma_0,C_0>0$, independent of $\nu$, such that
the following holds for all sufficiently small $\nu$.
If the perturbation satisfies \textup{(K3)} and
$\epsilon_{P_\nu}\le\gamma_0\alpha_\nu$, the KAM theorem produces a
Cantor set $\Pi_\nu^*\subset\Pi_\nu$ with
\begin{equation}\label{eq:uniform-measure}
 |\Pi_\nu\setminus\Pi_\nu^*|
 \le C_0(\operatorname{diam}\Pi_\nu)^{b-1}\alpha_\nu^\vartheta.
\end{equation}
\end{lemma}
\begin{proof}
Apply P\"oschel's Theorem~B first on the fixed reference box $\Pi^0$,
using the frequency bounds $M,L$.  This fixes a finite exceptional set
$\mathcal X_0$ and a constant $C_0$ in \eqref{eq:geometric}.  Under
restriction to $\Pi_\nu$, the domain-monotonicity statement in that
theorem gives an exceptional set $\mathcal X_\nu\subset\mathcal X_0$
and the same upper bound $C_0$.  The additional hyperbolic zones are
empty by Appendix~\ref{app:hyperbolic}.  It remains to treat the
elliptic resonance conditions indexed by $\mathcal X_\nu$.

Choose $\gamma_0$ according to Corollary~C so that the initial Fourier
cutoff satisfies $|k|\le K_0$ for every $(k,l)\in\mathcal X_0$.
This choice depends on the fixed set $\mathcal X_0$ and the fixed
frequency bounds.  Pairs in $\mathcal X_\nu$ then occur only at the
initial step, where the frequencies are those in \eqref{eq:frequencies}.
Their excluded zones are
\[
 \mathcal R^0_{k,l}(\alpha_\nu)
 =\left\{\xi\in\Pi_\nu:
     |f_{k,l}(\xi)|<\alpha_\nu
                    \frac{\langle l\rangle_1}{1+|k|^\tau}\right\},
 \qquad
 f_{k,l}(\xi)=k\cdot\omega(\xi)+l\cdot\Omega(\xi),
\]
where, in the increasing enumeration of the elliptic normal modes,
\[
 \langle l\rangle_1
   =\max\left(1,\left|\sum_{r\ge1}r\,l_{n_r}\right|\right).
\]

If $\mathcal X_0$ is nonempty, define
\[
 c_0=\frac12\min_{(k,l)\in\mathcal X_0}|f_{k,l}(0)|>0,\qquad
 D_0=\max_{(k,l)\in\mathcal X_0}\sup_{\xi\in\Pi^0}
                   \sum_{a=1}^b|\partial_{\xi_a}f_{k,l}(\xi)|<\infty.
\]
The positivity of $c_0$ follows from \eqref{eq:affine-nr}; the finiteness
of $D_0$ follows from the polynomial dependence on $\xi$ and the
finiteness of $\mathcal X_0$.  Choose $0<\nu_1\le\nu_*$ such that
$2D_0\nu_1\le c_0$.  For $\xi\in[0,2\nu_1]^b$,
\[
 |f_{k,l}(\xi)|\ge |f_{k,l}(0)|-D_0|\xi|_\infty\ge c_0,
 \qquad (k,l)\in\mathcal X_0.
\]
The weights $\langle l\rangle_1/(1+|k|^\tau)$ also have a finite
maximum on $\mathcal X_0$.  Since $\alpha_\nu\to0$, their products
with $\alpha_\nu$ are less than $c_0$ for sufficiently small $\nu$.
Every exceptional zone is then empty.  If $\mathcal X_0$ is empty,
this step is unnecessary.

Only the nonexceptional zones remain, and \eqref{eq:geometric} gives
\eqref{eq:uniform-measure}.  The frequency bounds, the lower bounds for the hyperbolic rates
and their pairwise differences, and the choice of $\gamma_0$
were fixed on $\Pi^0$, so the constants are
independent of $\nu$.
\end{proof}

\subsection{Choice of the normalization order}\label{sec:order-choice}
We determine a normalization order for which the perturbation estimates
satisfy the KAM smallness condition and the relative measure of the
excluded actions tends to zero.

If the normalization is stopped at an even degree $2m\ge4$, the same
finite induction gives a remainder $R_{2m+2}$ with
\[
 \|X_{R_{2m+2}}(w)\|_{3/2}\le C\|w\|_{1/2}^{2m+1}.
\]
On an action box $[\nu,2\nu]^b$, the tangential amplitude is
$O(\sqrt\nu)$.  Let $r$ be the normal radius, with $r^2\ll\nu$.
The terms of normal degree at least three give a scaled estimate
$O(\sqrt\nu\,r)$, whereas the action component of the remainder
gives $O(\nu^{m+1}/r^2)$.  The componentwise estimates below also give
the smaller terms, so the complete bound is
\begin{equation}\label{eq:general-order-estimate}
 \|X_P\|_r\le C\left(
 r^2+\sqrt\nu\,r+\frac{r^3}{\sqrt\nu}
       +\nu^m+\frac{\nu^{m+1}}{r^2}
                       +\frac{\nu^{m+1/2}}r\right).
\end{equation}
Balancing $\sqrt\nu\,r$ and $\nu^{m+1}/r^2$ gives
\begin{equation}\label{eq:general-order-balance}
 r=\nu^{(2m+1)/6},\qquad
 \|X_P\|_r\le C\nu^{(m+2)/3},\qquad
 \|X_P\|_r^{\Lip}\le C\nu^{(m-1)/3}.
\end{equation}
For $m\ge2$, this radius satisfies $r^2=o(\nu)$ and all the other
terms of \eqref{eq:general-order-estimate} have at least the same
power of $\nu$.  The Lipschitz estimate follows from one action
derivative, as verified below for $m=5$.

Put $\alpha=\nu^a$.  The analytic smallness condition requires
$a<(m+2)/3$, while \eqref{eq:geometric}, divided by $|\Pi_\nu|=\nu^b$,
gives a relative measure bound $C\nu^{a\vartheta-1}$.
The relative measure of the excluded actions therefore tends to zero if
\begin{equation}\label{eq:order-compatibility}
 \frac1\vartheta<a<\frac{m+2}{3}.
\end{equation}
Since $\vartheta=\eta/(1+\eta)$ approaches $1/2$ from below as
$\eta\uparrow1$, the smallest integer $m$ for which these estimates
are compatible with
\eqref{eq:order-compatibility} is $m=5$.  Thus the present perturbation
and measure estimates lead to normalization up to degree ten:
\[
 \begin{array}{c|cccc}
 \text{normalization degree }2m&4&6&8&10\\ \hline
 \text{remainder degree }2m+2&6&8&10&12\\
 \text{perturbation exponent }(m+2)/3
                              &4/3&5/3&2&7/3
 \end{array}
\]
For $m=5$, every $\beta<1/6$ is obtained by choosing $\eta$ close
enough to one and then $(1+\beta)/\vartheta<a<7/3$.

For comparison, the Schr\"odinger argument in
\cite{DongJianYuanNLS2026} uses quadratic spectral growth and a measure
bound linear in $\alpha$.  The corresponding inequalities are
$1<a<4/3$, which are already compatible at quartic order.

The higher resonant action terms are included in
\eqref{eq:frequencies} for the same reason.  For example, if the
degree-six tangential action polynomial were left in $P$, its angle
component would generally be $O(\nu^2)$.  Such a bound would not
satisfy $\|X_P\|_r=o(\nu^a)$ for the exponents $a>2$ required here.
Incorporating these terms changes the frequency map by $O(|\xi|^2)$;
\cref{prop:frequency-check} shows that this preserves the twist and
the required spectral bounds.

\subsection{Scaled perturbation estimates}
Let $\mathcal A_\xi$ denote \eqref{eq:AA} with the normal variables
unchanged.  Up to a scalar depending only on $\xi$, the transformed
Hamiltonian is $N_\xi+P_\nu$, where
\begin{equation}\label{eq:perturbation}
 \begin{split}
 P_\nu={}&Z_T(\xi+y)-Z_T(\xi)-\nabla Z_T(\xi)\cdot y\\
 &+\sum_{n\in\mathcal N_+}[Z_n(\xi+y)-Z_n(\xi)]z_n\widetilde z_n\\
 &+\sum_{h\in\Hyp}[Y_h(\xi+y)-Y_h(\xi)]a_hb_h
       +(P_{\ge3}+R_{12})\circ\mathcal A_\xi.
 \end{split}
\end{equation}
Choose
\begin{equation}\label{eq:radius}
 r_\nu=\nu^{11/6}.
\end{equation}
For $\nu$ small, $r_\nu^2<\nu/4$.  Hence the square roots in
\eqref{eq:AA} are holomorphic on $D(s,r_\nu)$, uniformly for
$\xi\in\Pi_\nu$.  Tangential amplitudes are $O(\sqrt\nu)$,
normal amplitudes are $O(r_\nu)$, and the total energy norm is
$O(\sqrt\nu)$.  As in \eqref{eq:scaled-norm}, normal vector fields
are measured in the $\cE_N^{3/2}$ norm.

\begin{lemma}\label{lem:perturbation-estimate}
On $D(s,r_\nu)\times\Pi_\nu$, $P_\nu$ satisfies \textup{(K3)} and
\begin{equation}\label{eq:P-estimates}
 \|X_{P_\nu}\|_{r_\nu}\le C\nu^{7/3},\qquad
 \|X_{P_\nu}\|_{r_\nu}^{\Lip}\le C\nu^{4/3}.
\end{equation}
\end{lemma}
\begin{proof}
Write $r=r_\nu$.  We estimate the three types of terms in
\eqref{eq:perturbation}.

\medskip
\noindent\textit{Step 1: Taylor remainders of the integrable part.}
Denote the first line of \eqref{eq:perturbation} by $T(\xi,y)$.
Taylor's formula gives
\[
 T(\xi,y)=\int_0^1(1-t)D^2Z_T(\xi+ty)[y,y]\dd t,
 \qquad
 \partial_yT=\int_0^1D^2Z_T(\xi+ty)y\dd t.
\]
The derivatives of $Z_T$ are bounded on a fixed complex neighborhood
of the reference box.  Hence $|\partial_yT|\le C|y|$, and the
only nonzero vector-field component of $T$ is the angle component.
Its scaled norm is therefore $O(r^2)$.  Differentiating the second
integral in $\xi$ gives the same bound for its parameter derivative.
For the elliptic mixed term, \eqref{eq:normal-coeff-bounds} implies
\[
 |Z_n(\xi+y)-Z_n(\xi)|\le C\langle n\rangle^{-1}|y|.
\]
Writing $d_n=Z_n(\xi+y)-Z_n(\xi)$, the corresponding normal
components are $-id_nz_n$ and $id_n\widetilde z_n$.  For example,
\[
 \sum_n\langle n\rangle^3|d_nz_n|^2
 \le C|y|^2\sum_n\langle n\rangle|z_n|^2.
\]
Thus this diagonal vector field maps $h^{1/2}$ to $h^{3/2}$
with norm at most $C|y|\|\zeta\|_{1/2}$, and its scaled norm is
$O(r^2)$.
The angle component is a convergent quadratic normal sum bounded by
$C\|\zeta\|_{1/2}^2$.  The finite hyperbolic mixed terms satisfy the
same estimates.  One parameter derivative uses only bounded higher
polynomial derivatives, so their Lipschitz seminorms are also $O(r^2)$.

\medskip
\noindent\textit{Step 2: terms of normal degree at least three.}
Let $t$ denote the tangential complex coordinates, let
$A=|t|+|\widetilde t|$, and set $R=\|\zeta\|_{1/2}$.
Bounded multilinearity, normal-degree decomposition, and the smoothing
of \cref{prop:NF} imply, when $A+R$ is sufficiently small,
\begin{equation}\label{eq:normal-degree-bounds}
 |\partial_tP_{\ge3}|+|\partial_{\widetilde t}P_{\ge3}|\le CR^3,
 \qquad
 \|X_{P_{\ge3}}^{\zeta}\|_{3/2}\le C(A R^2+R^3).
\end{equation}
For example, at degree four there is either one tangential and three
normal factors or four normal factors.  A normal derivative leaves
respectively one tangential and two normal factors or three normal
factors.  Higher total degrees introduce additional factors bounded by
powers of $A+R$ and satisfy the same estimates on a small ball.
The normal vector-field coefficients take values in $\cE_N^{3/2}$,
so the second estimate holds in that norm.

In action--angle coordinates,
$|\partial_y t|\le C\nu^{-1/2}$ and
$|\partial_\theta t|\le C\sqrt\nu$.
Taking $A\le C\sqrt\nu$ and $R\le r$ in
\eqref{eq:normal-degree-bounds} yields
\[
 |X^\theta|\le Cr^3\nu^{-1/2},\qquad
 r^{-2}|X^y|\le C\sqrt\nu\,r,\qquad
 r^{-1}\|X^\zeta\|_{3/2}\le C(\sqrt\nu\,r+r^2).
\]
For a tangential factor $t_a=(\xi_a+y_a)^{1/2}e^{\pm i\theta_a}$,
the parameter derivatives satisfy
\[
 |\partial_{\xi_a}t_a|\le C\nu^{-1/2},\qquad
 |\partial_{\xi_a}\partial_{y_a}t_a|\le C\nu^{-3/2},\qquad
 |\partial_{\xi_a}\partial_{\theta_a}t_a|\le C\nu^{-1/2}.
\]
Differentiating each of the finitely many tangential monomials therefore
introduces at most a factor $C\nu^{-1}$ in the preceding component
bounds.  Their normal multilinear coefficients are independent of
$\xi$.  Integrating the parameter derivatives along segments in the
convex box gives
\begin{equation}\label{eq:Pge3-scaled}
 \|X\|_r\le C(r^2+\sqrt\nu\,r+r^3/\sqrt\nu),\qquad
 \|X\|_r^{\Lip}\le C\nu^{-1}(r^2+\sqrt\nu\,r+r^3/\sqrt\nu).
\end{equation}

\medskip
\noindent\textit{Step 3: the remainder of order twelve.}
Let $V=X_{R_{12}}$ before the action--angle substitution.
By \eqref{eq:R12} and $\|w\|_{1/2}\le C\sqrt\nu$,
\[
 \|V\|_{3/2}\le C\nu^{11/2},\qquad
 \|DV\|_{\cL(\cE^{1/2},\cE^{3/2})}\le C\nu^5.
\]
For a tangential mode the exact conversion formulas are
\[
 X^{\theta_a}=\frac i2
       \left(\frac{V^{q_{j_a}}}{q_{j_a}}
                   -\frac{V^{\widetilde q_{j_a}}}{\widetilde q_{j_a}}\right),\qquad
 X^{y_a}=\widetilde q_{j_a}V^{q_{j_a}}+q_{j_a}V^{\widetilde q_{j_a}}.
\]
On a fixed complex strip, the tangential factors have size
$O(\sqrt\nu)$ and their reciprocals have size $O(\nu^{-1/2})$.
The normal components are unchanged.  Therefore
\begin{equation}\label{eq:R-scaled}
 \|X_{R_{12}\circ\mathcal A_\xi}\|_r
       \le C\left(\nu^5+\frac{\nu^6}{r^2}
                              +\frac{\nu^{11/2}}r\right).
\end{equation}
At fixed $(\theta,y,\zeta)$, $\|\partial_\xi w\|_{1/2}\le C\nu^{-1/2}$,
so $\|\partial_\xi V(w)\|_{3/2}\le C\nu^{9/2}$.
The parameter derivatives of a tangential amplitude and its
reciprocal are $O(\nu^{-1/2})$ and $O(\nu^{-3/2})$, respectively.
Differentiating the conversion formulas therefore gives
\[
 |\partial_\xi X^\theta|\le
 C(\nu^{-3/2}\nu^{11/2}+\nu^{-1/2}\nu^{9/2})\le C\nu^4,
 \qquad
 |\partial_\xi X^y|\le
 C(\nu^{-1/2}\nu^{11/2}+\nu^{1/2}\nu^{9/2})\le C\nu^5.
\]
The normal parameter derivative is bounded by $C\nu^{9/2}$ in
$\cE^{3/2}$.  The scaled Lipschitz estimate is consequently
\begin{equation}\label{eq:R-scaled-lip}
 \|X_{R_{12}\circ\mathcal A_\xi}\|_r^{\Lip}
       \le C\left(\nu^4+\frac{\nu^5}{r^2}
                              +\frac{\nu^{9/2}}r\right).
\end{equation}

For $r=\nu^{11/6}$ the largest terms in
\eqref{eq:Pge3-scaled} and \eqref{eq:R-scaled} are respectively
$\sqrt\nu\,r=\nu^{7/3}$ and $\nu^6/r^2=\nu^{7/3}$.
Their parameter bounds are $O(\nu^{4/3})$.
The Taylor remainders of the integrable part have size
$r^2=\nu^{11/3}$ and are smaller.
All maps are analytic on the common complex domain, and all normal
components lie in $h^{3/2}$ (with the finite hyperbolic coordinates
included).  This proves both \eqref{eq:P-estimates} and \textup{(K3)}.
\end{proof}

\subsection{Completion of the proof}
\begin{proof}[Proof of \cref{thm:main}]
Fix a nonlinearity $f$ satisfying \eqref{eq:nonlinearity}, a path in
\eqref{eq:universal-event}, an admissible $J$, and
$0<\beta<1/6$.  Since $(3/7)(1+\beta)<1/2$, choose $0<\eta<1$
sufficiently close to $1$ that
\[
 \vartheta=\frac\eta{1+\eta}>\frac37(1+\beta).
\]
Then choose an exponent $a$ with
\begin{equation}\label{eq:alpha-choice}
 \frac{1+\beta}{\vartheta}<a<\frac73,
 \qquad \alpha_\nu=\nu^a.
\end{equation}
Fix the strip width, the exponent \eqref{eq:tau}, the reference box,
and the bounds for the frequency maps.  For small $\nu$, the action--angle
image lies in the domain of the normal form and its remainder bounds.
By \cref{lem:perturbation-estimate},
\begin{equation}\label{eq:smallness-ratio}
 \frac{\epsilon_{P_\nu}}{\alpha_\nu}
 \le C\nu^{7/3-a}+\frac C M\nu^{4/3}\longrightarrow0.
\end{equation}
\Cref{prop:frequency-check,lem:uniform} give all frequency hypotheses
and a fixed positive smallness constant.  Thus the wave KAM theorem
applies for every sufficiently small $\nu$.
It gives a Cantor set $\Pi_\nu^*$ and a torus embedding
$\iota_\xi$ for each of its parameters.  The uniform measure estimate is
\[
 |\Pi_\nu\setminus\Pi_\nu^*|
 \le C\nu^{b-1}\alpha_\nu^\vartheta
 =C\nu^{b-1+a\vartheta}.
\]
Since $|\Pi_\nu|=\nu^b$ and $a\vartheta-1>\beta$,
this proves \eqref{eq:relative-measure} for $0<\nu<1$.
Taking the minimum of the finitely many smallness thresholds gives
$\nu_0(\omega,f,J,\beta)>0$.

\medskip
\noindent\textit{The tori in the original variables.}
The map from KAM coordinates to spectral coordinates is
$\Phi\circ\mathcal A_\xi$, because the normalized Hamiltonian is
$H\circ\Phi$.  It is a real analytic symplectic map, in local angle
charts.  The torus for the original Hamiltonian is therefore
\[
 \theta\longmapsto\Phi\circ\mathcal A_\xi\circ\iota_\xi(\theta),
\]
followed by the bounded linear map reconstructing $(u,v)$
from the spectral coordinates.
The KAM embedding estimate and \eqref{eq:smallness-ratio} imply an
$o(1)$ deformation of the trivial torus in scaled phase norm.
An angle displacement $o(1)$ changes a tangential spectral coordinate
by $o(\sqrt\nu)$; an action displacement $o(r_\nu^2)$ changes it by
$o(r_\nu^2/\sqrt\nu)=o(\sqrt\nu)$; and the normal displacement is
$o(r_\nu)=o(\sqrt\nu)$.
By \eqref{eq:Phi-bounds}, the Birkhoff map contributes only
$O(\nu^{3/2})$.  This proves the asserted energy-norm deformation of
\eqref{eq:linear-torus}.
After reducing $\nu_0$ if necessary, the displacement satisfies
$\sup_{\theta}\|u(\theta)\|_\infty<R_f/2$, so the nonlinearity is
defined along every trajectory on these tori.

All transformations preserve the reality conditions.
The normal variational equations in the KAM coordinates consist of
constant elliptic oscillators and finitely many hyperbolic pairs.  Their rates are small perturbations of the distinct
positive $\chi_h$, so their numbers of stable and unstable directions
are unchanged.  On the compact torus all conjugating differentials
and inverse differentials are bounded in the energy space.  They
therefore preserve the stated normal stability and hyperbolic splitting.

\medskip
\noindent\textit{Weak conjugacy and the wave equation.}
The trajectories on an analytic torus are $C^1$ in the energy space.
Let $\Psi=\Phi\circ\mathcal A_\xi$ and $K=H\circ\Psi$.  If $p(t)$
is such a trajectory in KAM coordinates and $q(t)=\Psi(p(t))$, then
symplecticity and the chain rule give, for every tangent vector $W$,
\[
 \Om(\dot q,W)
 =\Om(\dot p,D\Psi(p)^{-1}W)
 =\dd K(p)[D\Psi(p)^{-1}W]
 =\dd H(q)[W].
\]
This is the weak symplectic conjugacy used in
\cite{DongJianYuanNLS2026}.  In physical coordinates, testing the
resulting identity with arbitrary energy variations gives
\[
 u_t=v\quad\hbox{in }L^2,\qquad
 v_t=-A_\omega u-f(u)\quad\hbox{in }H^{-1}.
\]
Here $A_\omega:H^1_0\to H^{-1}$ is bounded by \eqref{eq:form}, and
$u\mapsto f(u)$ is analytic from a neighborhood of zero in $H^1_0$
to $H^1_0$, as proved in \cref{lem:smoothing}.
These identities prove the claimed time regularity and the PDE.

For completeness, define the spectral propagators
\[
 c_n(t)=\begin{cases}\cos(\mu_nt),&\lambda_n>0,\\
                     \cosh(\chi_nt),&\lambda_n<0,\end{cases}
 \qquad
 s_n(t)=\begin{cases}\sin(\mu_nt)/\mu_n,&\lambda_n>0,\\
                     \sinh(\chi_nt)/\chi_n,&\lambda_n<0.\end{cases}
\]
The finitely many negative modes cause only bounded growth on compact
time intervals.  The positive modes and the spectral energy norm
show that $c(t)$ preserves $H^1_0$ and $L^2$, while
$s(t):L^2\to H^1_0$ is bounded locally in time.
Variation of constants in each spectral coordinate gives
\begin{equation}\label{eq:mild}
 u(t)=c(t)u(0)+s(t)v(0)
              -\int_0^t s(t-s)f(u(s))\dd s.
\end{equation}
The integral is a Bochner integral in $H^1_0$.
The differentiated formula holds in $L^2$ and its derivative in
$H^{-1}$, using $s'=c$ and $c'=-A_\omega s$ in these spaces.
Thus \eqref{eq:mild} and the weak wave formulation are equivalent.
\end{proof}

\subsection{Analyticity in time}\label{sec:time-analytic}
We prove the additional regularity in \cref{thm:main}.
Let $\mathcal V_\xi:\T^b\to H^1_0\times L^2$ be the torus
embedding in the variables $(u,v)$.  Its analyticity and compactness of $\T^b$
give a bounded holomorphic extension to a complex strip
$|\Im\theta|_\infty<s_\xi$, after decreasing $s_\xi>0$ if necessary.
Write $\omega_*(\xi)\in\R^b$ for its rotation vector and put
\[
 \tau_\xi=\frac{s_\xi}{2(1+|\omega_*(\xi)|_\infty)}.
\]
For $|\Im t|<\tau_\xi$ the complexified trajectory
\[
 (u(t),v(t))=\mathcal V_\xi(\theta_0+\omega_*(\xi)t)
\]
is well defined and holomorphic with values in the complex energy
space, because its argument remains in $|\Im\theta|_\infty<s_\xi/2$.
Periodicity in the real angle variables makes this extension uniformly
bounded throughout the time strip.  If $M_\xi$ is such a bound, the
Banach-valued Cauchy formula on disks of radius $\tau_\xi/2$ gives
\[
 \sup_{t\in\R}\|\partial_t^k(u(t),v(t))\|_{H^1_0\times L^2}
       \le M_\xi k!\left(\frac2{\tau_\xi}\right)^k,
 \qquad k\ge0.
\]
On the real axis, the weak equation already proved identifies
$v=u_t$.  This proves the additional time analyticity asserted in
\cref{thm:main}.

\appendix
\section{KAM theorem and measure estimates for wave frequencies}\label{app:KAM}
We state the parts of P\"oschel's Theorems~A and~B and Corollary~C
\cite[pp.~120--123]{Poschel1996} used in Section~\ref{sec:KAM},
specialized to spectral exponent $d=1$ and zero exponential sequence
weight.  We then describe the extension to the finite hyperbolic modes
in \cite[Theorem~2.5]{YuanHyperbolic2003}.  The notation follows the
complex Hamiltonian coordinates of this paper.

\subsection{Hypotheses and existence theorem}
First suppose that all normal modes are elliptic.  On a compact
parameter set $\Pi\subset\R^b$ of positive measure, consider
\[
 H=N+P,\qquad
 N=\omega(\xi)\cdot y+\sum_{r\ge1}\Omega_r(\xi)z_r\widetilde z_r.
\]
Assume the following conditions.
\begin{enumerate}[label=\textup{(\Alph*)}]
\item The map $\omega$ is bi-Lipschitz onto its image.  For
$(k,l)\in\Z^b\times\Z^{(\N)}$, $(k,l)\ne0$, and $|l|\le2$,
the zero set of $k\cdot\omega+l\cdot\Omega$ has Lebesgue measure
zero.  Moreover, $l\cdot\Omega$ is nonzero throughout $\Pi$ for
$1\le|l|\le2$.
\item There is $0<\eta<1$ and a fixed sequence $\bar\Omega_r=r+c$
for all sufficiently large $r$, with arbitrary fixed initial entries,
such that $\widetilde\Omega=\Omega-\bar\Omega$ is Lipschitz into
the sequence space with norm
\[
 \|w\|_{\infty,\eta}=\sup_{r\ge1}r^\eta|w_r|.
\]
\item The perturbation is real analytic in the phase variables and
Lipschitz in $\xi$.  Its vector field maps the phase space with normal
index $p$ analytically into that with index $\bar p$, where
$\bar p-p\ge\eta>0$.  The supremum and parameter Lipschitz norms
of this vector field on $D(s,r)$ are finite.
\end{enumerate}
Here $D(s,r)$ and the scaled vector-field norm are defined as in
Section~\ref{sec:abstract-KAM}, with indices $p,\bar p$ in place of
$1/2,3/2$.  The Hamiltonian is real valued for real angles and actions
and $\widetilde z=\overline z$.  For a parameter-dependent map, the
superscript $\mathcal L$ denotes the sum of its supremum norm and
Lipschitz seminorm.  Choose $M,L\ge1$ such that
\[
 \|\omega\|_\Pi^{\mathcal L}
       +\|\widetilde\Omega\|_{\infty,\eta,\Pi}^{\mathcal L}\le M,
 \qquad
 \|\omega^{-1}\|_{\omega(\Pi)}^{\Lip}\le L.
\]
In P\"oschel's notation, $\delta=-\eta$ and
$p-\bar p\le\delta<0$.  We take $\tau>(b+3)(1+\eta)/\eta$,
which also satisfies the summability condition for the measure estimate.

\begin{theorem}[P\"oschel's Theorem A, wave case]\label{thm:app-A}
Under \textup{(A)}--\textup{(C)}, there is $\gamma>0$ such that
\[
 \epsilon_P=\|X_P\|_r+\frac\alpha M\|X_P\|_r^{\Lip}
       \le\gamma\alpha,\qquad 0<\alpha\le1,
\]
implies the existence of a Cantor set $\Pi_\alpha\subset\Pi$ and
a Lipschitz family of real-analytic invariant $b$-tori.  Their embeddings
$\iota_\xi$ extend holomorphically to $|\Im\theta|<s/2$ and carry
the motions $\iota_\xi(\theta+\omega_*(\xi)t)$, where $\omega_*$
is Diophantine.  The normal variational equations are reducible to
constant elliptic oscillators.  For
$\iota_0(\theta)=(\theta,0,0,0)$,
\begin{align}
 \|\iota-\iota_0\|_r
       +\frac\alpha M\|\iota-\iota_0\|_r^{\Lip}
 &\le C\epsilon_P/\alpha,\label{eq:app-embedding}\\
 \|\omega_*-\omega\|_{\Pi_\alpha}
       +\frac\alpha M\|\omega_*-\omega\|_{\Pi_\alpha}^{\Lip}
 &\le C\epsilon_P.\label{eq:app-frequencies}
\end{align}
The embedding norms are taken in the scaled phase norm on this strip.

The iteration also gives Lipschitz frequency maps
$\omega_m,\Omega_m$ on $\Pi$, $m\ge0$, with
$\omega_0=\omega$, $\Omega_0=\Omega$, and
\begin{align}
 \|\omega_m-\omega\|_\Pi
       +\frac\alpha M\|\omega_m-\omega\|_\Pi^{\Lip}
 &\le C\epsilon_P,\label{eq:app-tangential-correction}\\
 \|\Omega_m-\Omega\|_{\infty,\eta,\Pi}
       +\frac\alpha M\|\Omega_m-\Omega\|_{\infty,\eta,\Pi}^{\Lip}
 &\le C\epsilon_P.\label{eq:app-normal-correction}
\end{align}
Set
\[
 \begin{split}
 \mathcal Z&=\{(k,l)\in\Z^b\times\Z^{(\N)}:
                         (k,l)\ne0,\ |l|\le2\},\\
 \langle l\rangle_1&=\max\left(1,\left|\sum_{r\ge1}r l_r\right|\right),
 \qquad A_k=1+|k|^\tau,\\
 \mathcal Z_0&=\mathcal Z,\qquad
 \mathcal Z_m=\{(k,l)\in\mathcal Z:|k|>K_0 2^{m-1}\}\quad(m\ge1),
 \end{split}
\]
where $K_0\ge1$ is the initial Fourier cutoff.  With
\begin{equation}\label{eq:app-zones}
 \mathcal R^m_{k,l}(\alpha)
 =\left\{\xi\in\Pi:
 |k\cdot\omega_m(\xi)+l\cdot\Omega_m(\xi)|
                <\alpha\frac{\langle l\rangle_1}{A_k}\right\},
\end{equation}
the excluded parameters satisfy
\[
 \Pi\setminus\Pi_\alpha
 \subset\bigcup_{m\ge0}\ \bigcup_{(k,l)\in\mathcal Z_m}
                       \mathcal R^m_{k,l}(\alpha).
\]
\end{theorem}

The constants are independent of $r$ and $\alpha$ once the dimensions,
regularity indices, strip width, Diophantine exponent, and frequency
bounds have been fixed.  The normal stability conclusion follows from
the final normal form; see Remark~1 and Section~4 of
\cite{Poschel1996}.

\subsection{Measure estimates}
For $d=1$, the measure estimate uses the additional gap asymptotics
\begin{equation}\label{eq:app-gap}
 \frac{\Omega_r(\xi)-\Omega_s(\xi)}{r-s}
       =1+O(s^{-\eta})\qquad(r>s),
\end{equation}
uniformly for large $s$ and for $\xi\in\Pi$.

\begin{theorem}[P\"oschel's Theorem B, wave case]\label{thm:app-B}
Assume \textup{(A)}, \textup{(B)}, and \eqref{eq:app-gap}.
Suppose that, uniformly in $m$,
\begin{align*}
 \|\omega_m-\omega\|_\Pi&\le\alpha,
 &\|\Omega_m-\Omega\|_{\infty,\eta,\Pi}&\le\alpha,\\
 \|\omega_m-\omega\|_\Pi^{\Lip}&\le(2L)^{-1},
 &\|\Omega_m-\Omega\|_{\infty,\eta,\Pi}^{\Lip}&\le(2L)^{-1}.
\end{align*}
There are a finite set $\mathcal X\subset\mathcal Z$ and constants
$C_{\rm geom},\alpha_0>0$ such that, for $0<\alpha\le\alpha_0$,
\begin{equation}\label{eq:app-measure}
 \left|\bigcup_{m\ge0}\ \bigcup_{(k,l)\in\mathcal Z_m\setminus\mathcal X}
                     \mathcal R^m_{k,l}(\alpha)\right|
 \le C_{\rm geom}(\operatorname{diam}\Pi)^{b-1}
                  \alpha^{\eta/(1+\eta)}.
\end{equation}
These constants and $\mathcal X$ are determined by the unperturbed
frequency maps and their bounds.  The exceptional set and the measure
constant do not increase when the parameter domain is restricted to a
closed subset.
\end{theorem}

\begin{corollary}[P\"oschel's Corollary C]\label{cor:app-C}
The constant $\gamma$ in \cref{thm:app-A} can be reduced, depending
also on $M,L$ and $\mathcal X$, so that the frequency corrections
satisfy the hypotheses of \cref{thm:app-B} and
$K_0\ge\max_{(k,l)\in\mathcal X}|k|$ when $\mathcal X$ is nonempty.
For sufficiently small $\alpha$,
\begin{equation}\label{eq:app-total-measure}
 |\Pi\setminus\Pi_\alpha|
 \le C_{\rm geom}(\operatorname{diam}\Pi)^{b-1}
                         \alpha^{\eta/(1+\eta)}
       +\sum_{(k,l)\in\mathcal X}|\mathcal R^0_{k,l}(\alpha)|.
\end{equation}
For fixed $\Pi$ and fixed frequency maps, the right-hand side tends
to zero as $\alpha\downarrow0$.
\end{corollary}

Indeed, the choice of $K_0$ excludes every exceptional pair from
$\mathcal Z_m$ for $m\ge1$.  Its contribution is therefore confined
to $\mathcal R^0_{k,l}$, which uses the unperturbed frequencies.
Condition \textup{(A)} implies that the measure of each such zone tends
to zero.  No affine dependence on the parameters is required here.

In Section~\ref{sec:KAM}, the frequencies are polynomials in the
actions.  They are defined on the reference box $\Pi^0$, and
\eqref{eq:affine-nr} implies
$k\cdot\omega(0)+l\cdot\Omega(0)\ne0$ for every exceptional pair.
Lemma~\ref{lem:uniform} consequently makes all the finitely many initial
zones empty on $\Pi_\nu$ for sufficiently small $\nu$.
Theorem~\ref{thm:app-B} then gives \eqref{eq:uniform-measure}, with
constants chosen on $\Pi^0$.  The normal indices in this appendix
correspond to the increasing enumeration
$\mathcal N_+=\{n_1<n_2<\cdots\}$ in the main text.
The regularity indices are $p=1/2$, $\bar p=3/2$, and
$\delta=-\eta$.

\subsection{Finite hyperbolic modes}\label{app:hyperbolic}
Yuan's analytic theorem \cite[Theorem~2.5]{YuanHyperbolic2003}
allows finitely many imaginary normal frequencies in addition to the
real elliptic frequencies.  For the Hamiltonian \eqref{eq:abstract-N},
the canonical substitution
\[
 z_h=a_h,\qquad \widetilde z_h=-ib_h
\]
gives $i\dd z_h\wedge\dd\widetilde z_h=\dd a_h\wedge\dd b_h$
and normal frequency $i\Lambda_h$.  The reality conditions
are $a_h,b_h\in\R$ for these modes and
$\widetilde z_n=\overline{z_n}$ for elliptic modes.
The homological equations commute with the corresponding antilinear
involution.  Their normalized solutions, flows, and limiting embeddings
therefore preserve these conditions.
The diagonal normal-form corrections are consequently real in the
elliptic variables and in the products $a_hb_h$.  Thus the tangential
and elliptic frequencies remain real at each KAM step, as do the
hyperbolic rates.  The elliptic frequency maps therefore satisfy the
reality hypothesis of P\"oschel's measure theorem.

In this application the hyperbolic rates are positive and pairwise
distinct.  On the fixed reference box, choose $g_*>0$ such that
\[
 \Lambda_h\ge g_*,\qquad
 |\Lambda_h-\Lambda_{h'}|\ge g_*\quad(h\ne h').
\]
The second condition is omitted if there is only one hyperbolic pair.
Reducing the perturbation threshold preserves these inequalities with
$g_*/2$ for the corrected rates.  The existence and embedding conclusions
of \cref{thm:app-A} then hold with the finite hyperbolic variables
included in the phase norms.  The normal variational equations reduce
to constant elliptic oscillators and the same number of hyperbolic pairs.

The extra modes do not change the measure exponent in
\eqref{eq:app-measure}.  A divisor involving one hyperbolic frequency,
or two hyperbolic frequencies whose imaginary parts do not cancel,
has imaginary part bounded away from zero by the lower bounds for
the rates and their pairwise differences.  When two hyperbolic
contributions cancel, the divisor reduces to $k\cdot\omega_m$.  For a mixed elliptic--hyperbolic
divisor, its real part is $k\cdot\omega_m\pm\Omega_{n,m}$.
If $n>C(1+|k|)$, this part has size comparable to $n$, whereas the
exclusion threshold is at most $C\alpha n/(1+|k|^\tau)$.
Thus a nonempty zone would require $n\le C(1+|k|)$ for small
$\alpha$.  In that range the threshold is at most $C\alpha$,
which is smaller than the lower bound for the hyperbolic rates.  Such zones are therefore
empty.  Divisors involving only the finitely many hyperbolic modes
are treated using the same lower bounds.

Only the elliptic resonance zones remain in the parameter exclusion.
The measure estimate in Section~\ref{sec:KAM} follows from
P\"oschel's Theorem~B and Corollary~C.  Yuan's analytic extension
provides the invariant tori and their hyperbolic normal directions.

\section*{Declaration of AI use}
The authors used ChatGPT (OpenAI) for language editing and assistance
with symbolic calculations.  All mathematical arguments, calculations,
and results were checked and verified by the authors, who take full
responsibility for the content of this paper.


\begin{thebibliography}{99}

\bibitem{Bambusi2003}
D. Bambusi,
Birkhoff normal form for some nonlinear PDEs,
\emph{Comm. Math. Phys.} \textbf{234} (2003), 253--285.

\bibitem{BambusiGrebert2006}
D. Bambusi and B. Gr\'ebert,
Birkhoff normal form for partial differential equations with tame modulus,
\emph{Duke Math. J.} \textbf{135} (2006), 507--567.

\bibitem{DongJianYuanNLS2026}
Y. Dong, W. Jian, and X. Yuan,
Almost-sure quenched KAM tori for cubic NLS with a spatial white-noise potential,
\href{https://arxiv.org/abs/2609.22808}{arXiv:2609.22808}, 2026.

\bibitem{DongJianYuanSpectral2026}
Y. Dong, W. Jian, and X. Yuan,
High-energy asymptotics for finite-interval Schr\"odinger operators
with Gaussian white-noise potential,
\href{https://arxiv.org/abs/2606.22426}{arXiv:2606.22426}, 2026.

\bibitem{DumazLabbe2024}
L. Dumaz and C. Labb\'e,
Anderson localization for the 1-d Schr\"odinger operator with white noise potential,
\emph{J. Funct. Anal.} \textbf{286} (2024), no.~1, Paper No.~110191.

\bibitem{FukushimaNakao1977}
M. Fukushima and S. Nakao,
On spectra of the Schr\"odinger operator with a white Gaussian noise potential,
\emph{Z. Wahrscheinlichkeitstheorie verw. Gebiete} \textbf{37} (1977), 267--274.

\bibitem{GengSunWangWave2026}
J. Geng, Y. Sun, and W.-M. Wang,
Quasi-periodic solutions to nonlinear random wave equations at fixed potential realizations,
\emph{Math. Z.} \textbf{313} (2026), Paper No.~13.
\href{https://doi.org/10.1007/s00209-026-04031-6}{doi:10.1007/s00209-026-04031-6}.

\bibitem{GengSunWangKAM2026}
J. Geng, Y. Sun, and W.-M. Wang,
KAM theorem for discrete nonlinear random equations at fixed potential realizations,
\emph{Potential Anal.} \textbf{65} (2026), Paper No.~5.
\href{https://doi.org/10.1007/s11118-026-10279-8}{doi:10.1007/s11118-026-10279-8}.

\bibitem{HrynivMykytyuk2006}
R.~O. Hryniv and Ya.~V. Mykytyuk,
Eigenvalue asymptotics for Sturm--Liouville operators with singular potentials,
\emph{J. Funct. Anal.} \textbf{238} (2006), 27--57.

\bibitem{HrynivMykytyuk2012}
R.~O. Hryniv and Ya.~V. Mykytyuk,
Self-adjointness of Schr\"odinger operators with singular potentials,
\emph{Methods Funct. Anal. Topology} \textbf{18} (2012), 152--159.

\bibitem{Poschel1996}
J. P\"oschel,
A KAM-theorem for some nonlinear partial differential equations,
\emph{Ann. Scuola Norm. Sup. Pisa Cl. Sci. (4)} \textbf{23} (1996), 119--148.
\href{https://www.numdam.org/item/ASNSP_1996_4_23_1_119_0/}{Full text}.

\bibitem{PoschelWave1996}
J. P\"oschel,
Quasi-periodic solutions for a nonlinear wave equation,
\emph{Comment. Math. Helv.} \textbf{71} (1996), 269--296.

\bibitem{SavchukShkalikov2003}
A.~M. Savchuk and A.~A. Shkalikov,
Sturm--Liouville operators with distribution potentials,
\emph{Trans. Moscow Math. Soc.} (2003), 143--192.

\bibitem{Wayne1990}
C.~E. Wayne,
Periodic and quasi-periodic solutions of nonlinear wave equations via KAM theory,
\emph{Comm. Math. Phys.} \textbf{127} (1990), 479--528.

\bibitem{YuanHyperbolic2003}
X. Yuan,
Invariant manifold of hyperbolic-elliptic type for nonlinear wave equation,
\emph{Int. J. Math. Math. Sci.} \textbf{2003} (2003), no.~18, 1111--1136.
\href{https://doi.org/10.1155/S0161171203207092}{doi:10.1155/S0161171203207092}.

\bibitem{YuanWave2006}
X. Yuan,
Quasi-periodic solutions of nonlinear wave equations with a prescribed potential,
\emph{Discrete Contin. Dyn. Syst.} \textbf{16} (2006), no.~3, 615--634.
\href{https://doi.org/10.3934/dcds.2006.16.615}{doi:10.3934/dcds.2006.16.615}.

\bibitem{YuanSurvey2014}
X. Yuan,
KAM theorems and open problems for infinite-dimensional Hamiltonian with short range,
\emph{Sci. China Math.} \textbf{57} (2014), no.~7, 1479--1486.
\href{https://doi.org/10.1007/s11425-014-4779-1}{doi:10.1007/s11425-014-4779-1}.

\end{thebibliography}
\end{document}